%% file: generators.tex
\documentclass[11pt,reqno,tbtags,draft]{amsart}
\usepackage{amssymb}
\usepackage{url}
\usepackage[square,numbers]{natbib}
\bibpunct[, ]{[}{]}{;}{n}{,}{,}

\title[Generating Pareto records]
{Generating Pareto records}

\newcommand\urladdrx[1]{{\urladdr{\def~{{\tiny$\sim$}}#1}}}
\author{James Allen Fill}
\address{Department of Applied Mathematics and Statistics,
The Johns Hopkins University,
3400 N.~Charles Street,
Baltimore, MD 21218-2682 USA}
\email{jimfill@jhu.edu}
\urladdrx{http://www.ams.jhu.edu/~fill/}
\thanks{Research for both authors supported by
the Acheson~J.~Duncan Fund for the Advancement of Research in
Statistics.}

\author{Daniel~Q.\ Naiman}
\address{Department of Applied Mathematics and Statistics,
The Johns Hopkins University,
3400 N.~Charles Street,
Baltimore, MD 21218-2682 USA}
\email{daniel.naiman@jhu.edu}
\urladdrx{https://www.ams.jhu.edu/~dan/}

\numberwithin{equation}{section}

\renewcommand\le{\leqslant}
\renewcommand\ge{\geqslant}

\allowdisplaybreaks

\theoremstyle{plain}
\newtheorem{theorem}{Theorem}[section]
\newtheorem{lemma}[theorem]{Lemma}

\newtheorem{corollary}[theorem]{Corollary}
\newtheorem{conj}[theorem]{Conjecture}

\theoremstyle{definition}
\newtheorem{example}[theorem]{Example}
\newtheorem{definition}[theorem]{Definition}

\newtheorem{remark}[theorem]{Remark}

\newtheorem*{acks}{Acknowledgements}

\theoremstyle{remark}

\newenvironment{romenumerate}[1][-10pt]{
\addtolength{\leftmargini}{#1}\begin{enumerate}
 \renewcommand{\labelenumi}{\textup{(\roman{enumi})}}%
 \renewcommand{\theenumi}{\textup{(\roman{enumi})}}%
 }{\end{enumerate}}

\newcounter{oldenumi}
{\setcounter{oldenumi}{\value{enumi}}
\begin{romenumerate} \setcounter{enumi}{\value{oldenumi}}}
{\end{romenumerate}}

\newcounter{thmenumerate}

\newcounter{xenumerate}   

\newcommand\xfootnote[1]{\unskip\footnote{#1}$ $} 

\newcommand\pfitem[1]{\par(#1):}
\newcommand\pfitemx[1]{\par#1:}
\newcommand\pfitemref[1]{\pfitemx{\ref{#1}}}

\newcommand\pfcase[2]{\smallskip\noindent\emph{Case #1: #2} \noindent}
\newcommand\step[2]{\smallskip\noindent\emph{Step #1: #2} \noindent}
\newcommand\stepx{\smallskip\noindent\refstepcounter{steps}%
 \emph{Step \arabic{steps}:}\noindent}

\newcommand{\refT}[1]{Theorem~\ref{#1}}
\newcommand{\refC}[1]{Corollary~\ref{#1}}
\newcommand{\refL}[1]{Lemma~\ref{#1}}
\newcommand{\refR}[1]{Remark~\ref{#1}}
\newcommand{\refS}[1]{Section~\ref{#1}}
\newcommand{\refSS}[1]{Subsection~\ref{#1}}
\newcommand{\refP}[1]{Proposition~\ref{#1}}
\newcommand{\refD}[1]{Definition~\ref{#1}}
\newcommand{\refE}[1]{Example~\ref{#1}}
\newcommand{\refF}[1]{Figure~\ref{#1}}
\newcommand{\refApp}[1]{Appendix~\ref{#1}}

\newcommand{\refTab}[1]{Table~\ref{#1}}
\newcommand{\refand}[2]{\ref{#1} and~\ref{#2}}

\newcommand\marginal[1]{\marginpar{\raggedright\parindent=0pt\tiny #1}}
\newcommand\DN{\marginal{DN}}
\newcommand\JF{\marginal{JF}}

\newcommand\kolla{\marginal{CHECK! SJ}}
\newcommand\ms[1]{\texttt{[ms #1]}}
\newcommand\XXX{XXX \marginal{XXX}}
\newcommand\REM[1]{{\raggedright\texttt{[#1]}\par\marginal{XXX}}}
\newcommand\rem[1]{{\texttt{[#1]}\marginal{XXX}}}

\newcommand\linebreakx{\unskip\marginal{$\backslash$linebreak}\linebreak}

\begingroup
  \divide\count255 by 60
  \multiply\count255 by -60
  \advance\count255 by \time
  \ifnum \count255 < 10 \xdef\klockan{\the\count1.0\the\count255}
  \else\xdef\klockan{\the\count1.\the\count255}\fi
\endgroup

\newcommand\nopf{\qed}   
\newcommand\noqed{\renewcommand{\qed}{}} 
\newcommand\qedtag{\eqno{\qed}}

\DeclareMathOperator*{\sumx}{\sum\nolimits^{*}}
\DeclareMathOperator*{\sumxx}{\sum\nolimits^{**}}
\newcommand{\sumio}{\sum_{i=0}^\infty}
\newcommand{\sumjo}{\sum_{j=0}^\infty}
\newcommand{\sumj}{\sum_{j=1}^\infty}
\newcommand{\sumko}{\sum_{k=0}^\infty}
\newcommand{\sumk}{\sum_{k=1}^\infty}
\newcommand{\summo}{\sum_{m=0}^\infty}
\newcommand{\sumno}{\sum_{n=0}^\infty}
\newcommand{\sumn}{\sum_{n=1}^\infty}
\newcommand{\sumin}{\sum_{i=1}^n}
\newcommand{\sumjn}{\sum_{j=1}^n}
\newcommand{\sumkn}{\sum_{k=1}^n}
\newcommand{\prodin}{\prod_{i=1}^n}
\newcommand{\sprod}{\mbox{$\prod$}}

\newcommand\set[1]{\ensuremath{\{#1\}}}
\newcommand\bigset[1]{\ensuremath{\bigl\{#1\bigr\}}}
\newcommand\Bigset[1]{\ensuremath{\Bigl\{#1\Bigr\}}}
\newcommand\biggset[1]{\ensuremath{\biggl\{#1\biggr\}}}
\newcommand\lrset[1]{\ensuremath{\left\{#1\right\}}}
\newcommand\xpar[1]{(#1)}
\newcommand\bigpar[1]{\bigl(#1\bigr)}
\newcommand\Bigpar[1]{\Bigl(#1\Bigr)}
\newcommand\biggpar[1]{\biggl(#1\biggr)}
\newcommand\lrpar[1]{\left(#1\right)}
\newcommand\bigsqpar[1]{\bigl[#1\bigr]}
\newcommand\Bigsqpar[1]{\Bigl[#1\Bigr]}
\newcommand\biggsqpar[1]{\biggl[#1\biggr]}
\newcommand\lrsqpar[1]{\left[#1\right]}
\newcommand\xcpar[1]{\{#1\}}
\newcommand\bigcpar[1]{\bigl\{#1\bigr\}}
\newcommand\Bigcpar[1]{\Bigl\{#1\Bigr\}}
\newcommand\biggcpar[1]{\biggl\{#1\biggr\}}
\newcommand\lrcpar[1]{\left\{#1\right\}}
\newcommand\abs[1]{|#1|}
\newcommand\bigabs[1]{\bigl|#1\bigr|}
\newcommand\Bigabs[1]{\Bigl|#1\Bigr|}
\newcommand\biggabs[1]{\biggl|#1\biggr|}
\newcommand\lrabs[1]{\left|#1\right|}
\def\rompar(#1){\textup(#1\textup)}    
\newcommand\xfrac[2]{#1/#2}
\newcommand\xpfrac[2]{(#1)/#2}
\newcommand\xqfrac[2]{#1/(#2)}
\newcommand\xpqfrac[2]{(#1)/(#2)}
\newcommand\parfrac[2]{\lrpar{\frac{#1}{#2}}}
\newcommand\bigparfrac[2]{\bigpar{\frac{#1}{#2}}}
\newcommand\Bigparfrac[2]{\Bigpar{\frac{#1}{#2}}}
\newcommand\biggparfrac[2]{\biggpar{\frac{#1}{#2}}}
\newcommand\xparfrac[2]{\xpar{\xfrac{#1}{#2}}}
\newcommand\innprod[1]{\langle#1\rangle}
\newcommand\expbig[1]{\exp\bigl(#1\bigr)}
\newcommand\expBig[1]{\exp\Bigl(#1\Bigr)}
\newcommand\explr[1]{\exp\left(#1\right)}
\newcommand\expQ[1]{e^{#1}}
\def\xexp(#1){e^{#1}}
\newcommand\ceil[1]{\lceil#1\rceil}
\newcommand\floor[1]{\lfloor#1\rfloor}
\newcommand\lrfloor[1]{\left\lfloor#1\right\rfloor}
\newcommand\frax[1]{\{#1\}}
\newcommand\setn{\set{1,\dots,n}}
\newcommand\nn{[n]}
\newcommand\ntoo{\ensuremath{{n\to\infty}}}
\newcommand\Ntoo{\ensuremath{{N\to\infty}}}
\newcommand\asntoo{\text{as }\ntoo}
\newcommand\ktoo{\ensuremath{{k\to\infty}}}
\newcommand\mtoo{\ensuremath{{m\to\infty}}}
\newcommand\stoo{\ensuremath{{s\to\infty}}}
\newcommand\ttoo{\ensuremath{{t\to\infty}}}
\newcommand\xtoo{\ensuremath{{x\to\infty}}}
\newcommand\bmin{\wedge}
\newcommand\norm[1]{\|#1\|}
\newcommand\bignorm[1]{\bigl\|#1\bigr\|}
\newcommand\Bignorm[1]{\Bigl\|#1\Bigr\|}
\newcommand\downto{\searrow}
\newcommand\upto{\nearrow}
\newcommand\half{\tfrac12}
\newcommand\thalf{\tfrac12}

\newcommand\punkt{.\spacefactor=1000}    
\newcommand\iid{i.i.d\punkt}
\newcommand\ie{i.e\punkt}
\newcommand\eg{e.g\punkt}
\newcommand\viz{viz\punkt}
\newcommand\cf{cf\punkt}
\newcommand{\as}{a.s\punkt}
\newcommand{\aex}{a.e\punkt}
\newcommand{\io}{i.o\punkt}
\renewcommand{\ae}{\vu}  
\newcommand\whp{w.h.p\punkt}

\newcommand\ii{\mathrm{i}}

\newcommand{\tend}{\longrightarrow}
\newcommand\dto{\overset{\mathrm{d}}{\tend}}
\newcommand\pto{\overset{\mathrm{p}}{\tend}}
\newcommand\Pto{\overset{\mathrm{P}}{\tend}}
\newcommand\Lcto{\overset{\mathcal{L}}{\tend}}
\newcommand\asto{\overset{\mathrm{a.s.}}{\tend}}
\newcommand\eqd{\overset{\mathrm{d}}{=}}
\newcommand\neqd{\overset{\mathrm{d}}{\neq}}
\newcommand\op{o_{\mathrm p}}
\newcommand\Op{O_{\mathrm p}}

\newcommand\bbR{\mathbb R}
\newcommand\bbC{\mathbb C}
\newcommand\bbN{\mathbb N}
\newcommand\bbT{\mathbb T}
\newcommand\bbQ{\mathbb Q}
\newcommand\bbZ{\mathbb Z}
\newcommand\bbZleo{\mathbb Z_{\le0}}
\newcommand\bbZgeo{\mathbb Z_{\ge0}}

\newcounter{CC}
\newcommand{\CC}{\stepcounter{CC}\CCx} 
\newcommand{\CCx}{C_{\arabic{CC}}}     
\newcommand{\CCdef}[1]{\xdef#1{\CCx}}     
\newcommand{\CCname}[1]{\CC\CCdef{#1}}    
\newcommand{\CCreset}{\setcounter{CC}0} 
\newcounter{cc}
\newcommand{\cc}{\stepcounter{cc}\ccx} 
\newcommand{\ccx}{c_{\arabic{cc}}}     
\newcommand{\ccdef}[1]{\xdef#1{\ccx}}     
\newcommand{\ccname}[1]{\cc\ccdef{#1}}    
\newcommand{\ccreset}{\setcounter{cc}0} 

\renewcommand\Re{\operatorname{Re}}
\renewcommand\Im{\operatorname{Im}}

\newcommand\E{\operatorname{\mathbb E{}}}
\renewcommand\P{\operatorname{\mathbb P{}}}
\renewcommand\L{\operatorname{L}}
\newcommand\Var{\operatorname{Var}}
\newcommand\Cov{\operatorname{Cov}}
\newcommand\Corr{\operatorname{Corr}}
\newcommand\Exp{\operatorname{Exp}}
\newcommand\Po{\operatorname{Po}}
\newcommand\Bi{\operatorname{Bi}}
\newcommand\Bin{\operatorname{Bin}}
\newcommand\Be{\operatorname{Be}}
\newcommand\Ge{\operatorname{Ge}}
\newcommand\NBi{\operatorname{NegBin}}
\newcommand\Res{\operatorname{Res}}
\newcommand\fall[1]{^{\underline{#1}}}
\newcommand\rise[1]{^{\overline{#1}}}

\newcommand\supp{\operatorname{supp}}
\newcommand\sgn{\operatorname{sgn}}
\newcommand\Tr{\operatorname{Tr}}

\newcommand\ga{\alpha}
\newcommand\gb{\beta}
\newcommand\gd{\delta}
\newcommand\gD{\Delta}
\newcommand\gf{\varphi}
\newcommand\gam{\gamma}
\newcommand\gG{\Gamma}
\newcommand\gk{\varkappa}
\newcommand\gl{\lambda}
\newcommand\gL{\Lambda}
\newcommand\go{\omega}
\newcommand\gO{\Omega}
\newcommand\gs{\sigma}
\newcommand\gss{\sigma^2}
\newcommand\gth{\theta}
\newcommand\eps{\varepsilon}
\newcommand\ep{\varepsilon}

\newcommand\bJ{\bar J}

\newcommand\cA{\mathcal A}
\newcommand\cB{\mathcal B}
\newcommand\cC{\mathcal C}
\newcommand\cD{\mathcal D}
\newcommand\cE{\mathcal E}
\newcommand\cF{\mathcal F}
\newcommand\cG{\mathcal G}
\newcommand\cH{\mathcal H}
\newcommand\cI{\mathcal I}
\newcommand\cJ{\mathcal J}
\newcommand\cK{\mathcal K}
\newcommand\cL{{\mathcal L}}
\newcommand\cM{\mathcal M}
\newcommand\cN{\mathcal N}
\newcommand\cO{\mathcal O}
\newcommand\cP{\mathcal P}
\newcommand\cQ{\mathcal Q}
\newcommand\cR{{\mathcal R}}
\newcommand\cS{{\mathcal S}}
\newcommand\cT{{\mathcal T}}
\newcommand\cU{{\mathcal U}}
\newcommand\cV{\mathcal V}
\newcommand\cW{\mathcal W}
\newcommand\cX{{\mathcal X}}
\newcommand\cY{{\mathcal Y}}
\newcommand\cZ{{\mathcal Z}}

\newcommand\ett[1]{\boldsymbol1_{#1}}
\newcommand\bigett[1]{\boldsymbol1\bigcpar{#1}}
\newcommand\Bigett[1]{\boldsymbol1\Bigcpar{#1}}
\newcommand\etta{\boldsymbol1}

\newcommand\smatrixx[1]{\left(\begin{smallmatrix}#1\end{smallmatrix}\right)}

\newcommand\limn{\lim_{n\to\infty}}
\newcommand\limN{\lim_{N\to\infty}}

\newcommand\qw{^{-1}}
\newcommand\qww{^{-2}}
\newcommand\qq{^{1/2}}
\newcommand\qqw{^{-1/2}}
\newcommand\qqq{^{1/3}}
\newcommand\qqqb{^{2/3}}
\newcommand\qqqw{^{-1/3}}
\newcommand\qqqbw{^{-2/3}}
\newcommand\qqqq{^{1/4}}
\newcommand\qqqqc{^{3/4}}
\newcommand\qqqqw{^{-1/4}}
\newcommand\qqqqcw{^{-3/4}}

\newcommand\intii{\int_{-1}^1}
\newcommand\intoi{\int_0^1}
\newcommand\intoo{\int_0^\infty}
\newcommand\intoooo{\int_{-\infty}^\infty}
\newcommand\oi{[0,1]}
\newcommand\ooo{[0,\infty)}
\newcommand\ooox{[0,\infty]}
\newcommand\oooo{(-\infty,\infty)}
\newcommand\setoi{\set{0,1}}

\newcommand\dtv{d_{\mathrm{TV}}}

\newcommand\dd{\,\mathrm{d}}
\newcommand\ddx{\mathrm{d}}

\newcommand{\pgf}{probability generating function}
\newcommand{\mgf}{moment generating function}
\newcommand{\chf}{characteristic function}
\newcommand{\ui}{uniformly integrable}
\newcommand\rv{random variable}
\newcommand\lhs{left-hand side}
\newcommand\rhs{right-hand side}

\newcommand\gnp{\ensuremath{G(n,p)}}
\newcommand\gnm{\ensuremath{G(n,m)}}
\newcommand\gnd{\ensuremath{G(n,d)}}
\newcommand\gnx[1]{\ensuremath{G(n,#1)}}

\newcommand\etto{\bigpar{1+o(1)}}
\newcommand\GW{Galton--Watson}
\newcommand\GWt{\GW{} tree}
\newcommand\cGWt{conditioned \GW{} tree}
\newcommand\GWp{\GW{} process}
\newcommand\tX{{\widetilde X}}
\newcommand\tY{{\widetilde Y}}
\newcommand\kk{\varkappa}
\newcommand\spann[1]{\operatorname{span}(#1)}
\newcommand\tn{\cT_n}
\newcommand\tnv{\cT_{n,v}}
\newcommand\rea{\Re\ga}
\newcommand\wgay{{-\ga-\frac12}}
\newcommand\qgay{{\ga+\frac12}}
\newcommand\ex{\mathbf e}
\newcommand\xx{\mathbf x}
\newcommand\yy{\mathbf y}
\newcommand\zz{\mathbf z}
\newcommand\uu{\mathbf u}
\newcommand\hX{\hat X}
\newcommand\sgt{simply generated tree}
\newcommand\sgrt{simply generated random tree}
\newcommand\hh[1]{d(#1)}
\newcommand\WW{\widehat W}
\newcommand\coi{C\oi}
\newcommand\out{\gd^+}
\newcommand\zne{Z_{n,\eps}}
\newcommand\ze{Z_{\eps}}
\newcommand\gatoo{\ensuremath{\ga\to\infty}}
\newcommand\rtoo{\ensuremath{r\to\infty}}
\newcommand\Yoo{Y_\infty}
\newcommand\bes{R}
\newcommand\tex{\tilde{\ex}}
\newcommand\tbes{\tilde{\bes}}
\newcommand\Woo{W_\infty}
\newcommand{\hm}{m_1}
\newcommand{\thm}{\tilde m_1}
\newcommand{\bbb}{B^{(3)}}
\newcommand{\rr}{r^{1/2}}
\newcommand\coo{C[0,\infty)}
\newcommand\coT{\ensuremath{C[0,T]}}
\newcommand\expx[1]{e^{#1}}
\newcommand\gdtau{\gD\tau}
\newcommand\ygam{Y_{(\gam)}}
\newcommand\EE{V}
\newcommand\pigsqq{\sqrt{2\pi\gss}}
\newcommand\pigsqqw{\frac{1}{\sqrt{2\pi\gss}}}
\newcommand\gapigsqqw{\frac{(\ga-\frac12)\qw}{\sqrt{2\pi\gss}}}
\newcommand\gdd{\frac{\gd}{2}}

\newcommand\raisetagbase{\raisetag{\baselineskip}}

\newcommand\eit{e^{\ii t}}
\newcommand\emit{e^{-\ii t}}
\newcommand\tgf{\tilde\gf}
\newcommand\txi{\tilde\xi}
\newcommand\intT{\frac{1}{2\pi}\int_{-\pi}^\pi}
\newcommand\intpi{\int_{-\pi}^\pi}
\newcommand\Li{\operatorname{Li}}
\newcommand\doi{D_{01}}
\newcommand\cHoi{\cH(\doi)}
\newcommand\db{D}
\newcommand\dbm{D_-}
\newcommand\dbmb{\overline D_-}
\newcommand\dbmx{\widehat D_-}
\newcommand\bdb{\overline{D_B}}
\newcommand\xq{\setminus\set{\frac12}}
\newcommand\xo{\setminus\set{0}}
\newcommand\tqn{t/\sqrt n}
\newcommand\intpm[1]{\int_{-#1}^{#1}}
\newcommand\gnaxt{g_n(\ga,x,t)}
\newcommand\gaxt{g(\ga,x,t)}
\newcommand\gssx{\frac{\gss}2}
\newcommand\tq{\tilde q}
\newcommand\gao{\ga_0}
\newcommand\ppp{\cP_1}
\newcommand\dx{D^*}
\newcommand\tpsi{\tilde\psi}
\newcommand\tgD{\tilde\gD}
\newcommand\xinn{\xi_{n-1,N}}
\newcommand\zzn{\frac12+\ii y_n}
\newcommand\OHD{O_{\cH(D)}}
\newcommand\OHDx{O_{\cH(\dx)}}
\newcommand\tgdn{\tgD_N}
\newcommand\xgdn{\gD^*_N}
\newcommand\intt{\int_0^T}
\newcommand\act{|\cT|}

\newcommand{\ignore}[1]{}

\newcommand{\Holder}{H\"older}
\newcommand\CS{Cauchy--Schwarz}
\newcommand\CSineq{\CS{} inequality}
\newcommand{\Levy}{L\'evy}
\newcommand{\Takacs}{Tak\'acs}
\newcommand{\Frechet}{Fr\'echet}

\newcommand{\maple}{\texttt{Maple}}

\newcommand\citex{\REM}
\newcommand\refx[1]{\texttt{[#1]}}
\newcommand\xref[1]{\texttt{(#1)}}

\usepackage{tikz}
\usepackage{amssymb}
\usetikzlibrary{arrows,positioning}
\usepackage[english]{babel}
\usepackage{pgfplotstable}
\usepackage{pgfplots}
\usepackage{adjustbox}
\usetikzlibrary{calc}
\usepackage{relsize}
\usetikzlibrary{arrows.meta}
\tikzset{>={Latex[width=5mm,length=5mm]}}

\pgfplotsset{compat=1.3}

\usepackage{chngcntr}
\counterwithout{table}{section}

\begin{document}

\date{January~15, 2019}

\maketitle

\begin{abstract}
We present, (partially) analyze, and apply an efficient algorithm for the simulation of multivariate Pareto records.  A key role is played by minima of the record-setting region (we call these \emph{generators}) each time a new record is generated, and two highlights of our work are (i)~efficient dynamic maintenance of the set of generators and (ii)~asymptotic analysis of the expected number of generators at each time.
\end{abstract}

\section{Introduction}
\label{S:intro}

Records for univariate observations have been thoroughly studied; an excellent reference is the book by Arnold~\cite{Arnold(1998)}.  On the other hand, exploration of the fundamental behavior of multivariate records is relatively in its infancy, and simulations can provide highly useful empirical information.  The question thus arises: How does one sample multivariate records efficiently?  In this paper we provide and crudely analyze an algorithm based on importance sampling.  The algorithm's gain in efficiency from the naive approach of generating multivariate observations and waiting for (Pareto) records to occur is nothing short of dramatic: In the bivariate case, for example, it would take roughly $10^{61}$ observations to generate the 10,000 records (which result from just 10,000 fairly rapid steps of our algorithm) leading
to~\cite[Fig.~3]{Fillboundary(2018)} (and ultimately to the theoretical results of that paper) and roughly $10^{194}$ observations to generate the 100,000 records summarized
in \refTab{Table1} (which led to the still-incomplete theory in~\cite{Fillbreaking(2018)}).

We begin with notation and definitions, some repeated from~\cite[Sec.~1.1]{Fillboundary(2018)}.
For a positive integer~$n$, let $[n] := \{1, \dots, n\}$.  Thus $[d]^{[n]}$ denotes the set of all functions
from~$[n]$ into~$[d]$, or simply the set of all $n$-tuples with each entry in $\{1, \dots, d\}$.
For $d$-dimensional vectors $x = (x_1, \dots, x_d)$ and $y = (y_1, \dots, y_d)$,
write $x \prec y$ (respectively, $x \leq y$) to mean that $x_j < y_j$ (resp.,\ $x_j \leq y_j$)
for $j \in [d]$.
The notation $x \succ y$ means $y \prec x$, and $x \geq y$ means $y \leq x$; the notation $x < y$ means
$x \leq y$ but $x \neq y$, and $y > x$ means $x < y$.  We write
$x_+ := \sum_{j = 1}^d x_j$ (respectively, $\sprod x_j := \prod_{j = 1}^d x_j$) for
the sum (resp.,\ product) of coordinates of $x = (x_1, \dots, x_d)$.

We assume that underlying the creation of records are $d$-dimensional observations $X^{(1)}, X^{(2)}, \dots$ that are \iid\ (independent and identically distributed) copies of a random vector~$X$ with independent coordinates $X_j$.  Unless otherwise specified, we will assume that each $X_j$ has the uniform distribution over the interval $[0, 1)$.  Note
that it is simple (at least theoretically) to switch from uniform to any other distribution by means of the standard inverse probability transform (also known as the quantile transform).

\begin{definition}
\label{D:record}
(a)~We say that $X^{(k)}$ is a \emph{Pareto record} (or \emph{nondominated record} or \emph{weak record} or simply \emph{record}, or that $X^{(k)}$ \emph{sets} a record at time~$k$) if $X^{(k)} \not\prec X^{(i)}$ for all $1 \leq i < k$.

(b)~If $k \in [n]$,
we say that $X^{(k)}$ is a \emph{current record} (or \emph{remaining record}, or \emph{maximum}) at time~$n$ if $X^{(k)} \not\prec X^{(i)}$ for all $i \in [n]$.

(c)~If $k \in [n]$,
we say that $X^{(k)}$ is a \emph{broken record} at time~$n$ if it is a record but not a current record, that is, if $X^{(k)} \not\prec X^{(i)}$ for all $1 \leq i < k$ but $X^{(k)} \prec X^{(\ell)}$ for some $k < \ell \leq n$; in that case, the observation corresponding to the smallest such~$\ell$ is said to \emph{break} or \emph{kill} the
record $X^{(k)}$.
\end{definition}

\begin{definition}
\label{D:RS}
The \emph{record-setting region} at time~$n$ is the (random) closed set of points
\[
S_n := \{x \in [0, 1)^d: x \not\prec X^{(i)}\mbox{\ for all $i \in [n]$}\}.
\]
\end{definition}

Note that
\begin{align*}
S_n
&= \{x \in [0, 1)^d: x \not\prec X^{(i)}\mbox{\ for all $i \in [n]$} \\
&{} \qquad \qquad \qquad \mbox{such that $X^{(i)}$ is a current record at time~$n$}\},
\end{align*}
and that the current records at time~$n$ all belong to $S_n$ and lie on its
(topological)
boundary.

Our discussion of generators later in this paper will hinge on the notion of two specific kinds of orthants (translated and intersected with the unit hypercube $[0, 1)^d$).  For our purposes, it will be convenient to work in terms of \emph{closed} positive orthants and \emph{open} negative orthants (in the usual Euclidean topology, relativized to the hypercube).  Accordingly:

\begin{definition}
\label{D:orthant}
Suppose $x \in [0, 1)^d$.
\smallskip

\noindent
(a)~By the \emph{closed positive orthant generated} (or \emph{determined}) \emph{by~$x$}, we mean the set
\[
O^+_x := \{y \in [0, 1)^d: y \geq x\}.
\]

\noindent
(b)~By the \emph{open negative orthant generated} (or \emph{determined}) \emph{by~$x$}, we mean the set
\[
O^-_x := \{y \in [0, 1)^d: y \prec x\}.
\]
\end{definition}
\smallskip

To avoid large subscripts, we will sometimes use $O^+(x)$ [or $O^+(x_1, \dots, x_d)$] and $O^-(x)$ [or
$O^-(x_1, \dots, x_d)$] instead of $O^+_x$ and $O^-_x$, respectively.
Note that
\begin{equation}
\label{intersection}
O^+_x \cap O^+_y = O^+_{x \vee y},
\end{equation}
where $x \vee y$ is coordinatewise maximum: $x \vee y := (x_1 \vee y_1, \dots, x_d \vee y_d)$.
\smallskip

The contents of the paper are as follows.  We present an outline of our algorithm for sampling Pareto records, together with multivariate-records phenomena suggested from data collected by using the algorithm, in \refS{S:algo}.  A key part of the algorithm is the dynamic maintenance of the set of minima (the so-called ``generators'') of the record-setting region; discussion of this set maintenance begins in~\refS{S:generators}.  Good performance of the algorithm depends largely on \emph{efficient} maintenance of the set of generators; this issue is treated in \refS{S:efficiency}.  The complexity of the algorithm also depends on the number of generators each time a new record is generated; the counting of generators is discussed in Sections~\ref{S:char}--\ref{S:expected}.  More specifically, in \refS{S:char} we characterize generators in a way that facilitates their count; in \refS{S:deterministic} we derive a loose upper bound and a best-possible lower bound on the number of generators at any time in terms of the number of current records at that time; and in \refS{S:expected} we derive an exact expression (see Lemmas~\ref{L:GI} and~\ref{L:exact}) and an asymptotic expansion (\refT{T:Gdnasymptotics}) for the expected number of generators after the collection of any given number of observations.  (At the present, we are able to give only a fairly crude bound on the expected number  of generators in ``records-time'', that is, after the generation of any given number of Pareto records.)

\section{An algorithm for sampling Pareto records}
\label{S:algo}
In \refS{S:outline} we outline an algorithm for sampling Pareto records, and in \refS{S:empirical} we present two phenomena of records suggested from data we collected using the algorithm.

\subsection{An outline of the algorithm}
\label{S:outline}
The algorithm dynamically maintains a set~$G$ of ``generators'', such that the current record-setting region~$S$ is the union $\cup_{g \in G}\,O^+_g$.  (The elements of~$G$ are thus called generators because~$S$ is the up-set in $[0, 1)^d$ generated by~$G$ with respect to the partial order $\leq$.  Initially we have $G = \{0\}$.)

If we were generating observations one at a time, the next observation would set a record if and only if it were to fall in~$S$.  But we wish instead to generate a new Pareto record---or, equivalently, a new observation~$X$ (unconditionally distributed uniformly in $[0, 1)^d$) conditionally given that it falls in~$S$.  Note that the following two tasks are simple to carry out:
\begin{enumerate}
\item[(i)] For each $g \in G$, compute $\P(X \in O^+_g)$.  Indeed,
\[
\P(X \in O^+_g) = \sprod (1 - g_j).
\]
\item[(ii)] For any $g \in G$, sample~${\bf R}$ from the conditional distribution of~$X$ given $X \in O^+_g$.  Indeed, one need only sample a random vector~${\bf U}$ of independent Uniform$(0, 1)$ coordinates and set
\[
{\bf R}_j = g_j + (1 - g_j) {\bf U}_j, \quad j \in [d].
\]
\end{enumerate}

Now we are ready to outline our importance-sampling subroutine that takes as input the current set~$G$ of generators, outputs a new record~${\bf R}$, and updates~$G$ to $G'$ accordingly.  All sampling steps are to be performed mutually independently.
\bigskip

{\bf Subroutine for generating a new record~$R$:}
\begin{enumerate}
\item[1.] Sample~${\bf g}$ from~$G$ according to the distribution
\[
\P({\bf g} = g) = \frac{\P(X \in O^+_g)}{\sum_{h \in G} \P(X \in O^+_h)}
= \frac{\sprod (1 - g_j)}{\sum_{h \in G} \sprod (1 - h_j)}.
\]
\item[2.] Sample a random vector~${\bf U}$ of independent uniform$[0, 1)$ coordinates and set
\[
{\bf R}_j = {\bf g}_j + (1 - {\bf g}_j) {\bf U}_j, \quad j \in [d].
\]
\item[3.] If ${\bf R} = R$, accept~$R$ as a new record with probability
\[
1 / \#\{g: R \in O^+_g\}.
\]
If~$R$ is rejected, repeat Steps 1--3.
\smallskip
\item[4.] Update~$G$ to $G'$, as described in Sections~\ref{S:generators}--\ref{S:efficiency}.
\end{enumerate}
\ignore{
DN's slides (bivariate case), starting on PDF page 134.
\medskip

{\bf Goal:}~Sample from~$f$ conditional on being record-breaking pair
\medskip

page~143:
\smallskip

{\bf \large How to Sample from a Union:}
\medskip

{\bf Goal:}~To sample~$X$ according to pdf~$f$, conditioned on $X \in \cup_{i = 1}^n A_i$
\[
\tilde{f} := \frac{f I_{\cup_{i = 1}^n A_i}}{P_f\!\left( \cup_{i = 1}^n A_i \right)}
\]
where $P_f (B) = \int_B\!f(x)\,dx$
\medskip

{\bf Union Sampling Algorithm:}
\begin{itemize}
\item Sample~$I$ s.t.\ $P[I = i] = P_f(A_i) / \sum_j P_f (A_j)$
\item Sample $X \sim f$ conditioned on $X \in A_I$
\item Keep~$X$ with probability $1 / \#\{i : X \in A_i\}$
\end{itemize}
}
\bigskip

It is elementary to confirm that the ultimate output of Step~3 is indeed correctly distributed according to the conditional distribution of~$X$ given $X \in S$; and that (conditionally given~$S$, \ie,\ given~$G$) the probability of acceptance at each iteration of Step~3 is
\[
\frac{\P(X \in \cup_{g \in G} \,O^+_g)}{\sum_{g \in G} \P(X \in O^+_g)}
\geq \frac{\max_{g \in G} \P(X \in O^+_g)}{\sum_{g \in G} \P(X \in O^+_g)}
\geq \frac{1}{|G|},
\]
so that the number of times Steps 1--3 are repeated is distributed Geometric with success probability at least $1/|G|$ (and hence expected value at most $|G|$).

We will thus perform a conservative average-case analysis of Steps 1--3 of our subroutine in \refS{S:expected} (for which Sections \ref{S:char}--\ref{S:deterministic} are preparatory) by assessing the expected value of the number $|G|$ of generators at any given time.  The complexity of Step~4 will be discussed in Sections~\ref{S:generators}--\ref{S:efficiency}.

\begin{remark}
\label{R:importance}
The importance sampling algorithm builds on an algorithm for sampling from a union that appears in Karp and Luby~\cite{Karp(1983)} and Frigessi and Vercellis~\cite{Frigessi(1985)} and
is generalized in Naiman and Wynn~\cite{Naiman(1997)}, and whose extensions have been found to be useful for various applications, including genetics (Naiman and Wu~\cite{Naiman(2005)}) and medical imaging (Naiman and Priebe~\cite{Naiman(2001)}).
\end{remark}

\begin{figure}[htb]
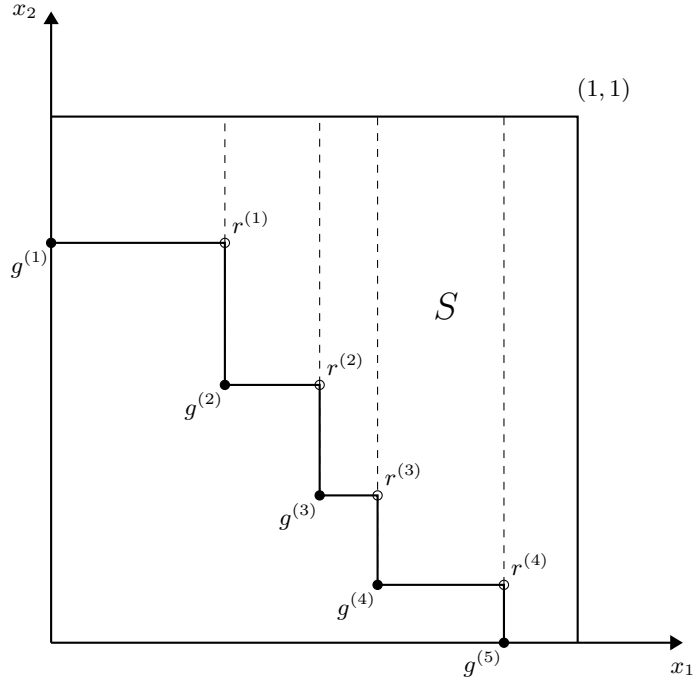

\include{figure0}
\caption{Example of a record-setting region~$S$ (above and to the right of the frontier formed from the current records and generators) in the bivariate case with four current records
($\rho=4$) and five generators ($\gamma=5$).
The record-setting region~$S$ is the union~\eqref{disjoint} of the rectangular regions marked with dashed lines.}
\label{F:bivariate}
\end{figure}

\begin{remark}
\label{R:bivariate}
An alternative to our importance-sampling algorithm is available in the bivariate case $d = 2$.  In that case we can write
$G = \{g^{(1)}, \ldots, g^{(\gamma)}\}$
with (\as)
\[
0 = g^{(1)}_1 < g^{(2)}_1 < \cdots < g^{(\gamma)}_1 < 1\quad\mbox{and}\quad
1 > g^{(1)}_2 > g^{(2)}_2 > \cdots > g^{(\gamma)}_2 = 0,
\]
and~$S$ is the \emph{disjoint} union
\begin{equation}
\label{disjoint}
S = \cup_{i = 1}^{\gamma} \left( [g^{(i)}_1, g^{(i + 1)}_1) \times [g^{(i)}_2, 1) \right)
\end{equation}
of rectangular regions, with the convention $g^{(\gamma + 1)}_1 = 1$.  Indeed, if we denote the set of current records by
$\{r^{(1)}, \ldots, r^{(\rho)}\}$
with
\[
0 < r^{(1)}_1 < r^{(2)}_1 < \cdots < r^{(\rho)}_1 < 1\quad\mbox{and}\quad
1 > r^{(1)}_2 > r^{(2)}_2 > \cdots > r^{(\rho)}_2 > 0,
\]
then $\gamma = \rho + 1$ and $g^{(i)} = (r^{(i - 1)}_1, r^{(i)}_2)$ for $i \in [\gamma]$ with the conventions $r^{(0)}_1 = 0$ and $r^{(\rho + 1)}_2 = 0$.  
See \refF{F:bivariate}.
With the representation~\eqref{disjoint} it is straightforward to sample uniformly from~$S$, and it is also rather simple to see (we omit the details) how to update~$G$ to $G'$.
\end{remark}

\subsection{Empirical phenomena}
\label{S:empirical}

We first noticed the phenomena of \refC{C:3d} and \refR{R:no formula} from data we collected using the algorithm outlined in the preceding subsection.  Moreover, two in-depth studies of multivariate-records behavior were also motivated empirically.

First, data such as those plotted
for $d = 2$ in~\cite[Fig.~3]{Fillboundary(2018)} suggested that, when observation coordinates have the standard exponential (rather than uniform) distribution, the topological boundary of the record-setting region after~$m$ records in general dimension $d \geq 2$ have been generated is likely to lie in a narrow strip between two parallel hyperplanes both located near the hyperplane $\{x: x_+ = (d! m)^{1/d}\}$.  This turned out to be provably true; see~\cite{Fillboundary(2018)}, especially
Thms.~5.3(a), 5.5(a), and 5.6(a) therein.

Second, long before this paper was written, the second-named author had used the algorithm to
produce \refTab{Table1}, with its striking Geometric($1/2$) pattern, for generation of 100,000 bivariate records.  A partial explanation of the pattern has now been provided in~\cite{Fillbreaking(2018)}.  Similarly, the following conjecture arises from data generated by the algorithm for larger dimensions:

\begin{center}
\begin{table}
\begin{tabular}{crc}
$k$&
$N_k$\ \ \ &
$\tilde{p}_k$ \\ \hline
0&50,334&0.50334\\
1&24,667&0.24667\\
2&12,507&0.12507\\
3&63,35&0.06335\\
4&3,040&0.03040\\
5&1,571&0.01571\\
6&782&0.00782\\
7&364&0.00364\\
8&202&0.00202\\
9&94&0.00094\\
10&48&0.00048\\
11&24&0.00024\\
12&18&0.00018\\
13&8&0.00008\\
14&4&0.00004\\
16&1&0.00001\\
17&0&0.00000\\
18&1&0.00001\\
\ & \ & \ \\
\end{tabular}
\caption{Results of simulation experiment in which $m =$\ 100,000 bivariate records are generated, and for each new record the number~$k$ of records it breaks is recorded.
The number of records that break~$k$ current records is denoted by $N_k$, and $\tilde{p}_k = N_k / m$ is the proportion of the 100,000 records that break~$k$ records.}
\label{Table1}
\end{table}
\end{center}

\begin{conj}
\label{Conj:kills}
Consider dimension $d \geq 2$.  Let $f_{d, m}$ denote the fraction of the first~$m$ records set that break~$0$ records.  Then there exist constants $p_d \in (0, 1)$ such that, almost surely, $f_{d, m} \to p_d$ as $m \to \infty$.  Further, $p_d \to 1$ as $d \to \infty$.
\end{conj}

The data also suggest that perhaps $p_d = 1 - d^{-1}$ for every $d \geq 2$; even for $d = 2$, the conjecture is stronger than what is proved in~\cite{Fillbreaking(2018)}.  For $d \geq 3$ and $k \geq 1$, we do not know what to conjecture concerning the limiting behavior of the fraction of the first~$m$ records set that break~$k$ records.

\section{The minima (aka generators) of the record-setting region}
\label{S:generators}

To sample records, all one needs to know (given the importance-sampling algorithm outlined in \refS{S:outline}) is how to express the record-setting region as a union of positive orthants (and how to update such an expression when a new record is generated)---and the more efficiently we do so, the more efficient the sampling will be.

Denote the set of current records at a given time by
$\{r^{(1)}, \dots, r^{(\rho)}\}$ (listed here in arbitrary, but fixed, order).
\smallskip

{\bf Assumption:\ }When discussing the deterministic geometry of such points $r^{(1)}, \dots, r^{(\rho)}$, we will for simplicity assume (as is almost surely true throughout time in our records model) that there are no ties in any coordinate, that is, that
\begin{equation}
\label{no coord ties}
\mbox{For each $j \in [d]$, the $\rho$ values $r^{(i)}_j$, $i \in [\rho]$, are distinct.}
\end{equation}
We will also assume that $r^{(i)} \in (0, 1)^d$ for $i \in [\rho]$, and that these points are (pairwise) incomparable; because we assume~\eqref{no coord ties}, for ``incomparable'' here we don't need to specify whether the partial order is $\preceq$ or $\leq$.
\smallskip

The record-setting region~$S$ is the closed set
\begin{equation}
\label{RS}
S = \cap_{i = 1}^{\rho} [O^-(r^{(i)})]^c.
\end{equation}
Algorithmically, what we want is to maintain the set
\[
G = \{g^{(1)}, \dots, g^{(\gamma)}\}
\]
of minima in the partial ordering $\leq$ of the record-setting region~$S$ (how we order~$G$ is irrelevant); we call these $g$'s the \emph{generators} of~$S$, because then
\[
S = \cup_{i = 1}^{\gamma} O^+_{g^{(i)}}
\]
is the minimal representation of~$S$ as a union of positive orthants.
So all we need to do is to describe (a)~how to initialize~$G$ and (b)~how to update~$G$ when a new record is set.
\refF{F:frontier_3d} shows an example of the \emph{frontier} (\ie,\ topological boundary) of the
record-setting region for $d=3$ when the number of current records is $\rho=8$ and
there are $\gamma=17$ generators.

For~(a), as noted previously in \refS{S:outline}, initially the record-breaking region is $S = O^+_0$ and has $\gamma = 1$ and $G = \{0\}$.

For (b), suppose that at the current time the record-setting region~$S$ has
$G = \{g^{(1)}, \dots, g^{(\gamma)}\}$ and a new record $r = (r_1, \dots, r_d)$ is generated.  Then
\[
S = \cup_{i = 1}^{\gamma} O^+_{g^{(i)}}
\]
and the updated~$S$, call it $S'$, is
\[
S' = S \cap (O^-_r)^c = \cup_{i = 1}^{\gamma} \left[ O^+_{g^{(i)}} \cap (O^-_r)^c \right],
\]
where complements are taken relative to the unit hypercube $[0, 1)^d = O^+_0$.  But
\begin{equation}
\label{key}
(O^-_r)^c = \cup_{k = 1}^d O^+_{r_k e^{(k)}},
\end{equation}
where, as usual, $e^{(k)} := (0, \dots, 0, 1, 0, \dots, 0)$ with the~$1$ as the $k$th coordinate.
Thus, for any $g \in [0, 1)^d$, using~\eqref{intersection} we find
\[
O^+_g \cap (O^-_r)^c
= \cup_{k = 1}^d \left( O^+_g \cap O^+_{r_k e^{(k)}} \right)
= \cup_{k = 1}^d O^+_{g \vee (r_k e^{(k)})}.
\]
Therefore,
 \[
S' = \cup_{i = 1}^{\gamma} \cup_{k = 1}^d O^+_{g^{(i)} \vee (r_k e^{(k)})}
\]
expresses the updated record-setting region as a union of closed positive orthants.  Observe, however, that the cardinality $d \gamma$ of the multiset
\[
\widehat{G} := \left\{ g^{(i)} \vee (r_k e^{(k)}): 1 \leq i \leq \gamma\mbox{\ and\ }1 \leq k \leq d \right\}
\]
is~$d$ times as large as the cardinality of the set~$G$.
To produce the updated set $G'$, one needs to find the minima of $\widehat{G}$ in the partial order $\leq$ on $O^+_0$.
\bigskip

\section{Improving efficiency}
\label{S:efficiency}

The algorithm we have described for finding generators can be made \emph{substantially} more efficient, as follows.  Say that a generator $g^{(i)}$ before the new record is generated \emph{survives} the new record~$r$ if $g^{(i)} \not\prec r$; which generators survive the new record can be determined simply by comparing each generator to the new record.  For each surviving $g^{(i)}$, we have $g^{(i)} \vee (r_j e^{(j)}) = g^{(i)}$ for some~$j$; thus we can
reduce~$\widehat{G}$ by replacing the~$d$ elements $g^{(i)} \vee (r_k e^{(k)})$ with $1\ \leq k \leq d$ by the single element $g^{(i)}$.  Further, $g^{(i)}$ does indeed ``survive'' as a minimum of $S'$, according to the following explanation.  After the new record is generated, $g^{(i)}$ still belongs to the record-setting region (because $g^{(i)} \not\prec r$), and it's still a minimum (because $S' \subset S$).  Further, the minima $g^{(i)}$ of~$S$ that are not surviving according to our definition do not ``survive'' (\ie, are ``killed''), in the sense that they don't even belong to~$S'$ [because $g^{(i)} \prec r$].

Call the set of surviving generators~$\Sigma$ and the set of non-surviving generators~$N$.  We have just now argued that $\Sigma \subseteq S'$.  As a next step, compare elements of the multiset
\[
\widehat{N} := \left\{ y \vee (r_k e^{(k)}): y \in N\mbox{\ and\ }1 \leq k \leq d \right\}
\]
in pairs to extract the set $N'$ of minima of $\widehat{N}$.  We claim that the disjoint union
\begin{equation}
\label{G'}
G' = \Sigma \cup N'
\end{equation}
is the desired set of generators of $S'$.

Here is a proof.  First we show that $N' \subseteq S'$, and for that consider $h = y \vee (r_k e^{(k)})$ with
$y \in N$.  Because $y \in N \subseteq S$, it follows for $1 \leq i \leq \gamma$ that $y \not\prec g^{(i)}$ and hence $h \not\prec g^{(i)}$.  Since it is also clear that $h \not\prec r$, it follows that $h \in S'$.

Now we need only show that $g < h$
(meaning that $g \leq h$ but $g \neq h$)
is not possible for $g \in \Sigma$ and
$h = y \vee (r_k e^{(k)}) \in N'$ (equivalently, for $h \in \widehat{N}$).  For that, first note that $g \not\leq y$, because~$g$ and~$y$ are distinct minima of~$S$.  So suppose $g_j > y_j$.  If we assume for the sake of contradiction that $g < y \vee (r_k e^{(k)})$ for some~$k$, then it must be that $k = j$.  But then $g_j \leq r_j$
and for every $\ell \neq j$ we have $g_{\ell} \leq y_{\ell} < r_{\ell}$, with the last inequality $y_{\ell} < r_{\ell}$ holding because $y \in N$.  Thus $g < r$.  If we can also argue that $g_j < r_j$, we will have the stronger conclusion that $g \prec r$, contradicting our assumption that $g \in \Sigma$.  To argue that $g_j < r_j$, observe [from the way that the set of generators evolves, confer~\eqref{r} below] that $g_j$ is the $j$th coordinate of some record other than~$r$.  Since $g_j \leq r_j$, we must have $g_j < r_j$, as desired.
\medskip

When what we have described is carried out in test dimension $d = 2$, it really does produce a decent gain in efficiency.  Let us explain why, in rough terms.  In dimension~$2$, if the number of remaining records killed by a new record is~$k$, then the number of elements of~$G$ killed by that record equals $k + 1$.  We know from~\cite{Fillbreaking(2018)} that (roughly put) with high probability the number of remaining records killed by a new record will be small (on the order of a constant).  Then, when we split~$G$ into the disjoint union of~$\Sigma$ and~$N$, the cardinality of~$N$ will with high probability be only on the order of a constant.  Thus if $|\Sigma| = s$ and $|N| = \nu$, the number of candidate elements of $S'$ will be only $s + d\,\nu = s + 2 \nu$ (only a constant more than~$g$, as opposed to cardinality $2 g$ for the multiset~$\widehat{G}$).  Moreover, the number of comparisons required to extract $G'$ \`a la~\eqref{G'} will be only
\[
\gamma + \mbox{${d \nu \choose 2}$} = \gamma + \mbox{${2 \nu \choose 2}$},
\]
where the first term counts the number of comparisons of old generators with the new record and the second term counts the comparisons required to extract $N'$ from~$\widehat{N}$.  Contrast this linear (in~$g$) count with the quadratic number of comparisons, namely, ${2 g \choose 2}$, required to compare every pair of elements in~$\widehat{G}$.

\section{Characterizing the generators of Pareto records}
\label{S:char}

Our development here begins by recalling~\eqref{RS}:
\begin{align}
S &= \cap_{i = 1}^{\rho} [O^-(r^{(i)})]^c
= \cap_{i = 1}^{\rho}  [\cup_{k = 1}^d O^+(r^{(i)}_k e^{(k)})] \nonumber \\
&= \cup_{k_1 = 1}^d \cdots \cup_{k_{\rho} = 1}^d \cap_{i = 1}^{\rho} O^+(r^{(i)}_{k_i} e^{(k_i)})
= \cup_{k_1 = 1}^d \cdots \cup_{k_{\rho} = 1}^d O^+(\vee_{i = 1}^{\rho} r^{(i)}_{k_i} e^{(k_i)}) \nonumber \\
\label{r}
&= \cup_{k \in [d]^{[\rho]}}\,O^+(R^{(\Pi_1(k))}_1, \ldots, R^{(\Pi_d(k))}_d),
\end{align}
where for $j \in [d]$ and $k \in [d]^{[\rho]}$ we have defined the ordered partition
$\Pi(k) = (\Pi_1(k), \dots, \Pi_d(k))$ of $[\rho]$ by
\[
\Pi_j(k) := k^{-1}(\{j\}) = \{i \in [\rho]:k_i = j\},
\]
and for $j \in [d]$ and $P \subseteq [\rho]$ we have defined
\[
R_j^{(P)} := \vee_{i \in P}\,r^{(i)}_j.
\]
Therefore we have the neat representation
\begin{equation}
\label{rep}
S = \bigcup O^+(R_1^{(\Pi_1)}, \ldots, R_d^{(\Pi_d)}),
\end{equation}
where the union here is taken over all ordered partitions $\Pi = (\Pi_1, \dots, \Pi_d)$ of $[\rho]$ into~$d$ sets; each $\Pi_j$  is allowed to be empty, in which case $R^{(\Pi_j)}_j := 0$.
This shows immediately that every element of~$G$ has in each coordinate either~$0$ or the value of some record in that coordinate.

To simplify our characterization of generators, we begin by considering only ``interior'' generators.
For any point $x \in O^+_0$, let $\tau(x)$ denote the set of non-zero coordinates of~$x$, and observe that~$x$ lies in the interior of $O^+_0$ if and only if $\tau(x) = [d]$.  We call such a point~$x$ an \emph{interior} point.

Observe that a point~$x$ of the form $(R_1^{(\Pi_1)}, \ldots, R_d^{(\Pi_d)})$ appearing in~\eqref{rep} is interior if and only if all the cells $\Pi$ of the partition are nonempty.  Next, note that $x \in (0, 1)^d$ is of such a form if and only if there exist~$d$ distinct indices $i_1, \dots, i_d$ such that $x_j = r^{(i_j)}_j$ for $j \in [d]$.

We are now in position to state and prove a characterization of the set~$I$ of interior generators.  (Note that
$I \subset G \subset S$.)

\begin{theorem}
\label{T:interior}
A point $g \in [0, 1)^d$ belongs to~$I$ if and only if
\begin{enumerate}
\item[(i)] $g \in S$, and
\item[(ii)] there exist~$d$ distinct indices $i_1, \dots, i_d$ such that
\begin{equation}
\label{min}
g_j = r^{(i_j)}_j = \min\{r^{(i_{\ell})}_j: \ell \in [d]\}\mbox{\ for every $j \in [d]$}.
\end{equation}
\end{enumerate}
\end{theorem}

\begin{proof}
First suppose $g \in I$.  Then~(i) is automatic from the definition of~$I$.  Moreover, we know from our earlier discussion that~(ii) holds for $g = (R_1^{(\Pi_1)}, \ldots, R_d^{(\Pi_d)})$ with the possible exception of the second equality in~\eqref{min}.  But if that equality does not hold, let $j, \ell \in [d]$ with $j \neq \ell$ satisfy
\[
r^{(i_{\ell})}_j < R_j^{(\Pi_j)}.
\]
We then move $i_{\ell}$ from the cell $\Pi_{\ell}$ to the cell $\Pi_j$ in order to form a new partition, call it $\Pi'$.  Then
\[
g > (R_1^{(\Pi'_1)}, \ldots, R_d^{(\Pi'_d)}) \in S,
\]
so~$g$ is not a generator.

Next we prove the converse.
If~$g$ has these two properties, then $g \in (0, 1)^d$ belongs to~$S$, so all that is left to show is that~$g$ is a minimum (with respect to~$\leq$) of~$S$.  Suppose that $x < g$;
we will complete the proof by showing that $x \not\in S$.

Let $j_0$ satisfy $x_{j_0} < g_{j_0}$.  Then
\begin{equation}
\label{j0}
x_{j_0} < g_{j_0} \leq r^{(i_{j_0})}_{j_0}
\end{equation}
using~\eqref{min} for the second inequality.  Additionally, for $j \neq j_0$ we have
\begin{equation}
\label{j}
x_j \leq g_j < r_j^{(i_{j_0})},
\end{equation}
where the second inequality holds by~\eqref{min} because
\[
g_j = r_j^{(i_j)} = \min\{r^{(i_{\ell})}_j: \ell \in [d]\},
\]
which is strictly smaller than $r_j^{(i_{j_0})}$ by our assumption~\eqref{no coord ties} because $i_j \neq i_{j_0}$.  Combining~\eqref{j0} and~\eqref{j}, we see that $x \prec r^{(i_{j_0})}$, and so $x \not\in S$.
\end{proof}

Now that we have characterized the interior generators, it is straightforward to characterize~$G$ in terms of projections of the current records to lower-dimensional coordinate subspaces, but some care must be taken to ensure that~\eqref{no coord ties} remains true after projection.  To begin a careful description, given a set $T = \{j_1, \dots, j_t\}$ of $[d]$ with $|T| = t \in [d]$ and $1 \leq j_1 < \cdots < j_t \leq d$, define the \emph{projection mapping} $\pi_T: \bbR^d \to \bbR^{t}$ by
\[
\pi_T(x_1, \dots, x_d) := (x_{j_1}, \dots, x_{j_t}),
\]
and define the \emph{injection mapping} $\iota_T: \bbR^t \to \bbR^d$ by
\[
\iota_T(x_1, \dots, x_t) = \vee_{k = 1}^t x_{j_k} e^{(j_k)}
\]
Recall that $\tau(x)$ denotes the set of nonzero coordinates of a point $x \in [0, 1)^d$.  Define the set of \emph{$T$-generators} to be the set
\[
G_T := G \cap \{x:\,\tau(x) = T\}
\]
and observe that~$G$ is the disjoint union
\[
G = \cup_{T \subseteq [d]} G_T.
\]
This observation, together with a characterization of each $G_T$, thus provides a characterization of~$G$.  A characterization of each $G_T$ is obtained by combining the following theorem with \refT{T:interior}.

To set up the statement of the theorem, consider the image
\[
R_T := \pi_T(R) = \{\pi_T(r^{(i)}): i \in [\rho]\} \subset \bbR^{|T|}
\]
under $\pi_T$ of the set $R := \{r^{(i)}: i \in [\rho]\}$ of current records, and note that $R_T$ inherits the property~\eqref{no coord ties} of ``no ties in any coordinate'' from~$R$.  Let $I_T$ denote the set of interior generators of $R_T$, and let $G'_T := \iota_T(I_T)$ denote the injection of $I_T$ into $\bbR^d$.

\begin{theorem}
\label{T:T}
For every $T \subseteq [d]$ we have $G_T = G'_T$.
\end{theorem}

\begin{proof}
As above, let $t = |T|$.
There is no loss of generality (and some ease in notation) in supposing that $|T| = [t]$, and thus $x \in G_T$ if and only if $x \in G$ and $x_{t + 1} = \cdots = x_d = 0$.
Let $x = (x_1, \dots, x_t, 0, \dots, 0)$ satisfy $\tau(x) = t$.  We will show that $x \in G_T$---equivalently, that
$x \in G$---if and only if $\pi_T(x) \in I_T$---equivalently, that $x \in \iota_T(I_T) = G'_T$.

Indeed, for~$x$ to be a generator, there are two requirements: (i)~$x \in S$, and (ii)~$x$ is a minimum of~$S$.  The requirement~(i) is that for each~$i$ there should exists $j \in [d]$ such that $x_j \geq r^{(i)}_j$.  However, since we assume that $r^{(i)} \succ 0$, such~$j$ must belong to $[t]$.  We have thus argued that
$x \in S = \mbox{RS}(R)$ if and only if $\pi_T(x) \in \mbox{RS}(R_T)$.

The requirement~(ii) is that $y < x$ must imply $y \notin S$.  But note that $y < x$ if and only if $y$ is of the form $y = (y_1, \dots, y_t, 0, \dots, 0)$ with $\pi_T(y) < \pi_T(x)$.  Thus requirement~(ii) can be rephrased thus:\ If $y = (y_1, \dots, y_t, 0, \dots, 0)$ with $\pi_T(y) < \pi_T(x)$, then $y \notin \mbox{RS}$---equivalently, by what we argued in connection with requirement~(i), that $\pi_T(y) \notin \mbox{RS}(R_T)$.

So we have argued that~$x$ is a generator if and only if $\pi_T(x) \in I_T$, \ie,\ if and only if $x \in G'_T$.  This is as desired.
\end{proof}

\begin{figure}[htb]
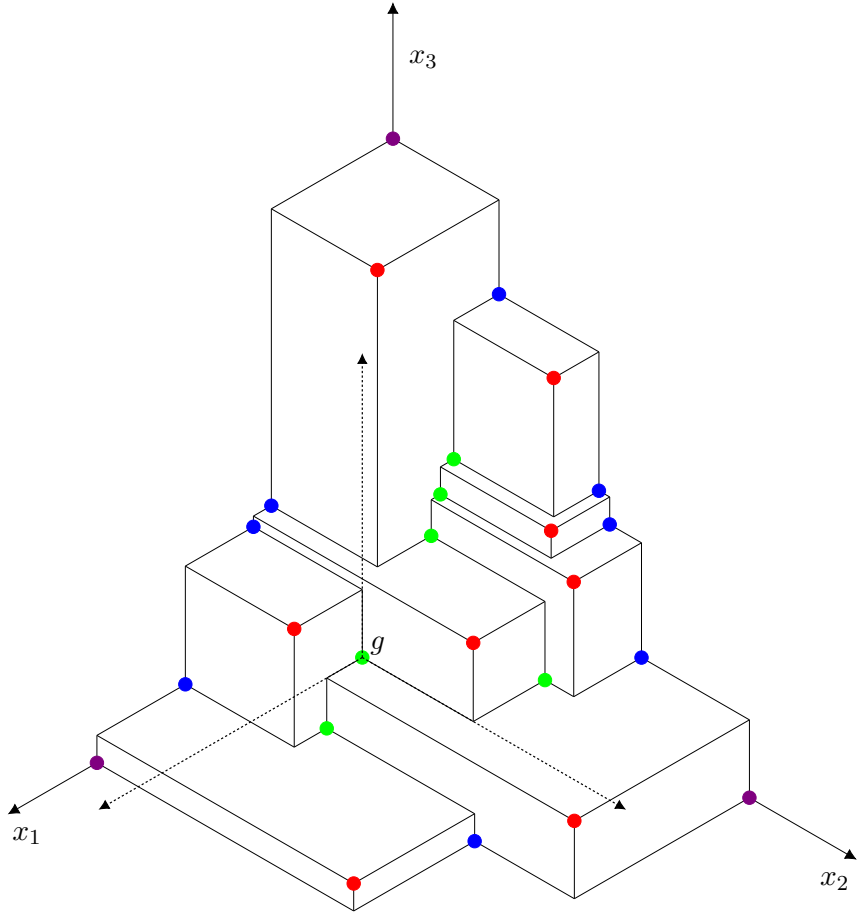

\include{figure1}
\caption{Example of a record frontier in dimension $d=3$ with $\rho=8$ remaining records shown in red and the resulting $\gamma = 17$ generators:\
three one-dimensional generators shown in violet, eight two-dimensional generators shown in blue, and six three-dimensional (interior) generators shown in green.
The lower boundary of one of the orthants $O^+_g$ is shown using dashed lines.}
\label{F:frontier_3d}
\end{figure}

In light of \refT{T:T}, we call the number of nonzero coordinates of a generator its \emph{dimension}.
\refF{F:frontier_3d} shows the generators
of various dimensions for an example with $d=3.$

\begin{example}
Suppose
$d = 4$ and the current records are $(0.2, 0.8, 0.3, 0.7)$ and $(0.5, 0.1, 0.4, 0.6)$.  Then $|G| = 8$, because $|G_T| = 1$ for precisely eight nonempty subsets~$T$ of $[4]$ and $|G_T| = 0$ otherwise.  The eight subsets~$T$ for which $|G_T| = 1$ are
\begin{align*}
G_{\{1\}} &= \{(0.5, 0, 0, 0)\};\,G_{\{2\}} = \{(0, 0.8, 0, 0)\};\,G_{\{3\}} = \{(0, 0, 0.4, 0)\}; \\
G_{\{4\}} &= \{(0, 0, 0, 0.7)\};\,
G_{\{1, 2\}} = \{(0.2, 0.1, 0, 0)\};\,G_{\{1, 4\}} = \{(0.2, 0, 0, 0.6)\}; \\
G_{\{2, 3\}} &= \{(0, 0.1, 0.3, 0)\};\,
G_{\{3, 4\}} = \{(0, 0, 0.3, 0.6)\}.
\end{align*}
Thus there are four one-dimensional generators, four two-dimensional generators, and no generators with dimension exceeding two.
\end{example}

\section{Deterministic upper and lower bounds \\ on the number of generators}
\label{S:deterministic}

We are interested, at each time unit, in bounding the number~$\gamma$ of (current) generators in terms of the number of the number~$\rho$ of current records, because we already know a fair amount about the latter in our records study.  In \refS{S:ub} we produce a crude upper bound.  In \refS{S:lb} we obtain a lower bound that in every dimension~$d$ is tight in the sense that there exists a set of~$\rho$ points that have all $d \rho$ coordinates distinct and are incomparable with respect to $\prec$ such that the number of generators achieves the lower bound.

\subsection{Upper bound}
\label{S:ub}

Currently, the best (almost sure) \emph{upper} bound we know on~$\gamma$ for general dimension~$d$ is the following.
\begin{theorem}
\label{T:ub}
We have
\begin{equation}
\label{ub}
\gamma \leq {{\rho + d - 1} \choose {d - 1}} \sim \frac{\rho^{d - 1}}{(d - 1)!},
\end{equation}
with asymptotics here for fixed~$d$ as $\rho \to \infty$.
\end{theorem}

\begin{proof}
We obtain this from the preceding \refS{S:char} [see especially~\eqref{rep}] as follows.
For a generator, we may suppose that~$\Pi_1$ contains the indices of all current records with the $\ell_1$ smallest first coordinates, for some $0 \leq \ell_1 \leq \rho$; this follows by the same argument as in the characterization of generators in the preceding section.  Having fixed $\ell_1$, we may now suppose for a generator that~$\Pi_2$ contains the current-record indices in $\Pi_1^c$ having the $\ell_2$ smallest second coordinates, for some $0 \leq \ell_2 \leq \rho - \ell_1$, and so on.  It follows that~$\gamma$ is at most the number of $d$-part compositions $\ell = (\ell_1, \dots, \ell_d)$ into~$d$ nonnegative parts.  That number is the binomial-coefficient bound in~\eqref{ub}.
\end{proof}

\begin{remark}
While the inequality in~\eqref{ub} is tight for $d = 1$ (then $\gamma = \rho = 1$) and for $d = 2$, it is quite loose for $d = 3$ (cf.~\refC{C:3d}).  It may well be for each $d \geq 4$ that there exists a universal constant $c_d$ such that $\gamma \leq c_d \rho$, regardless of the set of current records, but we do not know how to prove such a result.  However, for a result in this direction, see \refT{T:Gdnasymptotics} and parts~(a) and~(c) of \refR{R:time} following it.
\end{remark}

\subsection{Lower bound}
\label{S:lb}

\begin{theorem}
\label{T:lb}
For any dimension $d \geq 1$ and any given number $\rho \geq 0$ of current records (with all $d \rho$ coordinates distinct), the smallest number~$\gamma$ of generators possible is
\[
\gamma_0(\rho) := (d - 1) \rho + 1.
\]
\end{theorem}

\begin{remark}
\label{R:stronger}
Concerning achievement of the lower bound, the proof will establish the stronger result that there exists a sequence of incomparable points such that, for \emph{every} $\rho \geq 0$, the number $\gamma(\rho)$ of current generators after the first~$\rho$ points have arrived satisfies $\gamma(\rho) = \gamma_0(\rho)$.
\end{remark}

\begin{proof}[Proof of \refT{T:lb}]
(a)~We first show that
\begin{equation}
\label{lb}
\gamma \geq \gamma_0(\rho),
\end{equation}
and our proof is by induction on $\rho$.  For $\rho = 0$ (and for $\rho = 1$), this inequality (and indeed the conclusion of the theorem) is clear.  To complete the proof of~\eqref{lb} it therefore suffices to show that if a new $(\rho + 1)$st current record is added without killing any current records, then the number of generators increases by at least $d - 1$.

We use notation consistent with that of \refS{S:efficiency}.  For present generator-counting purposes, we may assume without loss of generality that the new record $r = r^{(\rho + 1)}$ doesn't kill any of the current records
$r^{(i)}, i \in [\rho]$, and that $r_1 < r_1^{(i)}$ for all $i \in [\rho]$.  Because of this, all interior generators in the set~$G$ of current generators survive~$r$; what's more, if a generator $g \in G$ is killed by~$r$, then
$g_1 = 0$.

Suppose now that exactly~$\nu$ ($\geq 1$) generators in~$G$ are killed by~$r$, and write
\[
N = \{g^{(1)}, \dots, g^{(\nu)}\}.
\]
Let
\[
\widehat{N}_k := \{g^{(i)} \vee (r_k e^{(k)}): i \in [\nu]\}\mbox{\ \ for $k \in [d]$}
\]
and
\[
\widehat{N} = \cup_{k = 1}^d \widehat{N}_k.
\]
The size of each set $\widehat{N}_k$ is obviously~$\nu$.  Moreover, the sets $\widehat{N}_k$ are \emph{pairwise incomparable sets}, in the sense that if $h_j \in \widehat{N}_j$ and $h_k \in \widehat{N}_k$ with $j \neq k$, then $h_j \not\leq h_k$ and $h_k \not\leq h_j$.  If $N'_k$ denotes the set of minima of
$\widehat{N}_k$, it therefore follows that the sets $N'_k$ are disjoint with union equal to the set $N'$ of minima of $\widehat{N}$.

Let $\nu'_k := |N'_k| \geq 1$. Let $\gamma' = |G'|$ count the number of post-new-record generators, and recall that the proof of~\eqref{lb} is complete if we can show $\gamma' \geq \gamma + d - 1$.  We have established that
\[
\gamma' = \gamma - \nu + \sum_{k = 1}^d \nu'_k.
\]
But a key observation is that
\[
\widehat{N}_1 = N = N'_1
\]
and hence $\nu'_1 = \nu$.  We conclude that
\begin{equation}
\label{gamma'}
\gamma' = \gamma + \sum_{k = 2}^d \nu'_k \geq \gamma + d - 1,
\end{equation}
completing the proof of~\eqref{lb}.

\begin{figure}[htb]
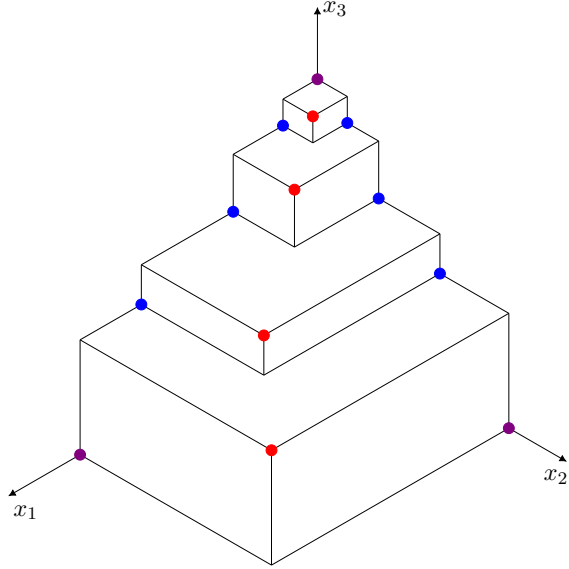

\include{figure2}
\caption{Example of $\rho = 4$ remaining records in dimension $d = 3$ achieving the lower bound $\gamma_0(\rho) = (d - 1) \rho + 1 = 9$ in \refT{T:lb}.
The four remaining records are shown in red, the three one-dimensional generators are shown in violet, and the six two-dimensional generators are shown in blue.}
\label{F:achieving}
\end{figure}

(b)~Here is an example as promised in \refR{R:stronger}.
(See \refF{F:achieving}.)
Suppose that the first coordinates of the successive records decrease strictly but all other coordinates increase strictly.  It is then easy to check [by induction, by following part~(a) of this proof] that the generators after~$\rho$ records have arrived are as follows.  There are precisely $d$ one-dimensional generators, namely, $r^{(1)}_1 e^{(1)}$ and $r^{(\rho)}_k e^{(k)}$ for $k = 2, \dots, d$.  There are precisely $(\rho - 1) (d - 1)$ two-dimensional generators, namely,
$r^{(i)}_1 e^{(1)} \vee r^{(i - 1)} e^{(k)}$ for $i = 2, \dots, \rho$ and $k = 2, \dots, d$.  There are no generators of higher dimension.  Thus the total number of generators is
\[
d + (d - 1)(\rho - 1) = (d - 1) \rho + 1 = \gamma_0(\rho),
\]
as desired.
\end{proof}

\begin{remark}
One might hope that the equality in~\eqref{gamma'} could also be used to improve the upper bound of \refT{T:ub}, but we don't know how to achieve this when $d \geq 4$.  (For $d = 3$, see \refC{C:3d} below.)
\end{remark}
\medskip

When $d = 1$, we always have $\gamma = \rho = 1$.  When $d = 2$, we always have $\gamma = \rho + 1$.  Our next result gives such a formula for $d = 3$.  The result is not actually a corollary of \refT{T:lb}, but its proof extends that of the theorem.

\begin{corollary}
\label{C:3d}
When $d = 3$, it is always true that
\[
\gamma = 2 \rho + 1.
\]
\end{corollary}

\begin{proof}
We return to the proof of \refT{T:lb}, recalling the equality in~\eqref{gamma'}.  Moreover, we see that
\[
N'_2 = \{(0, r_2, g^{(i_2)}_3)\} \quad \mbox{and} \quad N'_3 = \{(0, g^{(i_3)}_2, r_3)\}
\]
where $i_2$ is chosen to minimize $g^{(i)}_3$ over $i \in [\nu]$ and $i_3$ is chosen similarly to minimize
$g^{(i)}_2$.  The upshot is that updating by one record as in the proof of \refT{T:lb} we have
\[
\gamma' = \gamma + 2;
\]
That is, increasing~$\rho$ by~$1$ always increases~$\gamma$ by~$2$.  Since $\gamma = 1$ when
$\rho = 0$, the asserted result follows.
\end{proof}

\begin{remark}
\label{R:no formula}
There is no formula such as that of \refC{C:3d} (at least not for all values of~$\rho$) for any dimension
$d \geq 4$.  Indeed, it is easy to check that for $\rho = 2$ and $d \geq 1$, the distinct possible values
for~$\gamma$ are precisely those of the form
\[
\gamma = d + a (d - a)
\]
with $1 \leq a \leq \lfloor d / 2 \rfloor$.  In particular, the smallest possible value is $2 d - 1$, as guaranteed by \refT{T:lb}, and the largest possible value is $d + \lfloor d / 2 \rfloor \lceil d / 2 \rceil$.

As examples, when $d = 4$ there are two possible values, namely, $\gamma = 7$ and $\gamma = 8$; and when $d = 7$ there are three possible values, namely, $\gamma = 13$, $\gamma = 17$, and $\gamma = 19$.
\end{remark}

\section{The expected number of generators}
\label{S:expected}

In
this section, we abbreviate
$\dd x_1 \cdots \dd x_d$ as $\dd \xx$ and natural logarithm as~$\L$.

\subsection{Exact expressions}
\label{S:exact}

Let $\cG_{d, n}$ (respectively, $\cI_{d, n}$) denote the expected number of generators (resp.,\ interior generators) after a given number~$n$ of $d$-dimensional observations.  Our first result relates these two quantities.

\begin{lemma}
\label{L:GI}
For integers $d \geq 0$ and $n \geq 0$, we have
\begin{equation}
\label{Gdn}
\cG_{d, n} = \sum_{k = 0}^d {d \choose k} \cI_{k, n},
\end{equation}
\end{lemma}

\begin{proof}
This is immediate from \refT{T:T} and preceding discussion.
\end{proof}

In \refL{L:GI}, note that $\cI_{0, n} = \delta_{0, n}$: There is a single $0$-dimensional generator (namely, the origin in $\bbR^d$) if $n = 0$ and no $0$-dimensional generators otherwise.

The next result gives an exact expression for $\cI_{d, n}$.
We write $n^{\underline{k}}$ for the falling factorial power
\[
n (n - 1) \cdots (n - k + 1) = k! \mbox{${n \choose k}$}.
\]
\begin{lemma}
\label{L:exact}
For integers $d \geq 1$ and $n \geq 0$, we have
\begin{equation}
\label{Idn}
\cI_{d, n}
= n\fall{d} \int_{[0, 1)^d}\,\left[ \sprod (1 - x_j) \right]^{d - 1} \left[ 1 - \sprod (1 - x_j) \right]^{n - d} \dd \xx
\end{equation}
\end{lemma}

\begin{proof}
Referring to \refT{T:interior}(ii), let us say that the $d$-tuple $(X^{(i_1)}, \ldots, $ $X^{(i_d)})$ of observations (where the indices $i_j$ are distinct elements of $[n]$) \emph{generates} an interior generator~$g$ if
\[
g_j = X^{(i_j)}_j = \min\{X^{(i_{\ell})}_j: \ell \in [d]\}\mbox{\ for every $j \in d$}.
\]
Note that every interior generator is generated by precisely one such generating $d$-tuple.  Thus $\cI_{d, n}$ equals $n\fall{d}$ times the probability that $(X^{(1)}, \ldots, X^{(d)})$ generates an interior generator.  Condition on the value $\xx := (x^{(1)}, \ldots, x^{(d)})$ of this $d$-tuple.  According to \refT{T:interior}, in order for $\xx$ to generate an interior generator, two conditions are required.  One is that
\begin{equation}
\label{ellj}
x^{(\ell)}_j \geq x^{(j)}_j \mbox{\ for every $\ell, j \in d$ with $\ell \neq j$}.
\end{equation}
Let $x := (x^{(1)}_1, \ldots, x^{(d)}_d)$.
The other condition is that the remaining $n - d$ observations each need to fall outside $O^+_x$, guaranteeing the condition $x \in S$ required by \refT{T:interior}(i).

Therefore,
\[
\cI_{d, n}
= n\fall{d} \int_{\xx:\mbox{\scriptsize \eqref{ellj} holds}} \left[ 1 - \sprod (1 - x^{(j)}_j) \right]^{n - d}
\dd x^{(1)}_1 \dd x^{(1)}_2 \cdots \dd x^{(d)}_{d - 1} \dd x^{(d)}_d,
\]
a $d^2$-dimensional integral which reduces effortlessly to the $d$-dimensional integral~\eqref{Idn}.
\end{proof}

\subsection{Asymptotics}
\label{S:asymptotics}

From here we follow the same outline as for the expected number of remaining records in Bai et al.~\cite{Bai(2005)} to obtain an asymptotic expansion for $\cI_{d, n}$ (see our \refT{T:Gdnasymptotics}, the main result of \refS{S:expected}).  Accordingly, we begin by considering a Poissonized analogue of $\cI_{d, n}$.

\begin{lemma}
\label{L:Poissonized}
For integers $d \geq 1$ and $n \geq 0$, define
\[
\widehat{\cI}_{d, n}
:= n^d \int_{[0, 1)^d} \left[ \sprod (1 - x_j) \right]^{d - 1} \exp[- n\,\sprod (1 - x_j)] \dd \xx.
\]
Then, for fixed~$d$, as $n \to \infty$ we have
\[
\widehat{\cI}_{d, n}
= (\L n)^{d - 1} \sum_{j = 0}^{d - 1} \frac{(-1)^j \Gamma^{(j)}(d)}{j! (d - 1 - j)!} (\L n)^{-j}
+ O((n \L n)^{d - 1} e^{-n}).
\]
\end{lemma}

\begin{proof}
We have the following simple derivation:
\begin{align*}
\widehat{\cI}_{d, n}
&= n^d \int_{[0, 1)^d} \left( \sprod y_j \right)^{d - 1} \exp(- n\,\sprod y_j) \dd \yy
\mbox{\quad ($x_j \equiv 1 - y_j$)} \\
&= \int_{[0, n^{1/d})^d} \left( \sprod u_j \right)^{d - 1} \exp(- \sprod u_j) \dd \uu
\mbox{\quad ($y_j \equiv n^{-1/d} u_j$)} \\
&= \int_{(-d^{-1} \L n, \infty)^d} \exp \left[ - d\,z_+ - e^{- z_+} \right] \dd \zz
\mbox{\quad ($u_j \equiv e^{- z_j}$)} \\
&= \frac{1}{(d - 1)!} \int_{- \L n}^{\infty}\!(\L n + x)^{d - 1} \exp\left( - d\,x - e^{-x} \right) \dd x
\mbox{\quad ($x = z_+$)} \\
&= \frac{1}{(d - 1)!} \int_0^n\!(\L n - \L y)^{d - 1} y^{d - 1} e^{-y} \dd y
\mbox{\quad ($x = - \L y$)} \\
&= \frac{(\L n)^{d - 1}}{(d - 1)!} \sum_{j = 0}^{d - 1} {{d - 1} \choose j} \frac{(-1)^j}{(\L n)^j}
\int_0^n\!(\L y)^j y^{d - 1} e^{-y} \dd y \\
&= (\L n)^{d - 1} \sum_{j = 0}^{d - 1} \frac{(-1)^j}{j! (d - 1 - j)!} (\L n)^{-j} \int_0^n\!(\L y)^j y^{d - 1} e^{-y} \dd y.\\
\end{align*}
But
\begin{align*}
 \int_0^n\!(\L y)^j y^{d - 1} e^{-y} \dd y
&=  \int_0^{\infty}\!(\L y)^j y^{d - 1} e^{-y} \dd y - \int_n^{\infty}\!(\L y)^j y^{d - 1} e^{-y} \dd y \\
&=  \Gamma^{(j)}(d) - \Theta((\L n)^j n^{d - 1} e^{-n}).
\end{align*}
The desired result now follows easily.
\end{proof}

We next bound the difference between $\widehat{\cI}_{d, n}$ and
\begin{equation}
\label{Idntilde}
\widetilde{\cI}_{d, n}
:= n^d \int_{[0, 1)^d} \left[ \sprod (1 - x_j) \right]^{d - 1} \left[ 1 - \sprod (1 - x_j) \right]^n \dd \xx.
\end{equation}

\begin{lemma}
\label{L:quadratic}
For fixed $d \geq 1$, as $n \to \infty$ we have
\[
0 \leq \widehat{\cI}_{d, n} - \widetilde{\cI}_{d, n} = O(n^{-1} (\L n)^{d - 1}).
\]
\end{lemma}

\begin{proof}
We utilize the elementary inequality
\[
e^{- n t} (1 - n t^2) \leq (1 - t)^n \leq e^{- n t}
\]
for $n \geq 1$ and $0 \leq t \leq 1$ (see~\cite[Lemma~5]{Bai(2001)}).  This yields
\[
0 \leq \widehat{\cI}_{d, n} - \widetilde{\cI}_{d, n}
\leq n^{d + 1} \int_{[0, 1)^d} \left[ \sprod (1 - x_j) \right]^{d + 1} \exp[- n\,\sprod (1 - x_j)] \dd \xx.
\]
Proceeding just as in the proof of \refL{L:Poissonized}, we find that the last expression here is
$O(n^{-1} (\L n)^{d - 1})$.
\end{proof}

\begin{theorem}
\label{T:Idn}
For fixed $d \geq 1$, as $n \to \infty$ the expected number $\cI_{d, n}$ of interior generators at time~$n$ in dimension~$d$ satisfies
\[
\cI_{d, n}
= (\L n)^{d - 1} \sum_{j = 0}^{d - 1} \frac{(-1)^j \Gamma^{(j)}(d)}{j! (d - 1 - j)!} (\L n)^{-j} + O(n^{-1} (\L n)^{d - 1}).
\]
\end{theorem}

\begin{proof}
Comparing~\eqref{Idn} and~\eqref{Idntilde} and then invoking \refL{L:quadratic}, we see that
\begin{align*}
\cI_{d, n}
&= \frac{n\fall{d}}{(n - d)^d} \widetilde{\cI}_{d, n - d}
= [1 + O(n^{-1})]\,\widetilde{\cI}_{d, n - d} \\
&= [1 + O(n^{-1})]\,\left[ \widehat{\cI}_{d, n - d} + O(n^{-1} (\L n)^{d - 1}) \right] \\
&= [1 + O(n^{-1})]\,\widehat{\cI}_{d, n - d} + O(n^{-1} (\L n)^{d - 1}).
\end{align*}
But, according to \refL{L:Poissonized},
\begin{align*}
\widehat{\cI}_{d, n - d}
&= [\L (n - d)]^{d - 1} \sum_{j = 0}^{d - 1} \frac{(-1)^j \Gamma^{(j)}(d)}{j! (d - 1 - j)!} [\L (n - d)]^{-j}
+ O((n \L n)^{d - 1} e^{-n}) \\
& = (\L n)^{d - 1} \sum_{j = 0}^{d - 1} \frac{(-1)^j \Gamma^{(j)}(d)}{j! (d - 1 - j)!} (\L n)^{-j}
+ O(n^{-1} (\L n)^{d - 2}).
\end{align*}
Thus
\begin{align*}
\cI_{d, n}
&= [1 + O(n^{-1})]\,(\L n)^{d - 1} \sum_{j = 0}^{d - 1} \frac{(-1)^j \Gamma^{(j)}(d)}{j! (d - 1 - j)!} (\L n)^{-j}
+ O(n^{-1} (\L n)^{d - 1}) \\
&= (\L n)^{d - 1} \sum_{j = 0}^{d - 1} \frac{(-1)^j \Gamma^{(j)}(d)}{j! (d - 1 - j)!} (\L n)^{-j}
+ O(n^{-1} (\L n)^{d - 1}),
\end{align*}
as claimed.
\end{proof}

Combining~\eqref{Gdn} and~\eqref{Idn}, we can obtain an exact expression for $\cG_{d, n}$.  Similarly, combining~\eqref{Gdn} and \refT{T:Idn} we obtain the following asymptotic expansion in powers of logarithm for $\cG_{d, n}$ after a little rearrangement.

\begin{theorem}
\label{T:Gdnasymptotics}
For fixed $d \geq 1$, as $n \to \infty$ the expected number $\cG_{d, n}$ of generators at time~$n$ in dimension~$d$ satisfies
\[
\cG_{d, n}
= (\L n)^{d - 1} \sum_{j = 0}^{d - 1} a_{d, j} (\L n)^{-j} + O(n^{-1} (\L n)^{d - 1}),
\]
where
\[
a_{d, j}
:= \sum_{k = 0}^j {d \choose {d - j + k}} \frac{(-1)^{k} \Gamma^{(k)}(d - j + k)}{k! (d - 1 - j)!}.
\nopf
\]
\end{theorem}

\begin{remark}
\label{R:time}
(a)~In particular, $a_{d, 0} = 1$, so $\cG_{d, n}$ has lead-order asymptotics
\[
\cG_{d, n}
= (\L n)^{d - 1} + O((\L n)^{d - 2});
\]
this is $(d - 1)!$ times as large as the lead-order asymptotics for the expected number of remaining records, namely,
\[
\cR_{d, n}
= \frac{(\L n)^{d - 1}}{(d - 1)!} + O((\L n)^{d - 2}).
\]

(b)~We hope to extend the work of this section by finding at least lead-order asymptotics for the variance, and also a normal approximation or other limit theorem, for the number of generators after~$n$ observations.

(c)~\refT{T:Gdnasymptotics}
concerns the number of generators after~$n$ observations have been collected, but analysis of the algorithm presented in \refS{S:algo} requires us to know about the number of generators after~$m$ records have been generated.
It seems difficult to produce reasonably sharp theory concerning the latter number, call it $\tilde{\gamma}_m$. Denote the corresponding number of remaining records by $\tilde{\rho}_m$, and let $T_m$ denote the number of observations needed to set~$m$ records. One might conjecture
\begin{equation}
\label{conjecture}
\mbox{\as}:\ \tilde{\gamma}_m \sim (d - 1)! \tilde{\rho}_m \sim (d! m)^{(d - 1) / d},
\end{equation}
where the first asymptotic equivalence is suggested by part~(a) of this remark and the second
(known to be true for $d \geq 5$)
by conjecture~(5.3) in~\cite{Fillboundary(2018)}.  At present, the most we can assert with proof
for $d \geq 5$ is that
\begin{equation}
\label{as bound 5}
\mbox{\as}:\ \tilde{\gamma}_m
\leq (1 + o(1)) \frac{\tilde{\rho}_m^{d - 1}}{(d - 1)!}
= (1 + o(1)) \frac{(d! m)^{(d - 1)^2 / d}}{[(d - 1)!]^d},
\end{equation}
where the inequality follows from \refT{T:ub} and the equality from (5.3) in~\cite{Fillboundary(2018)}.
For $d = 4$, the situation is even slightly worse than~\eqref{as bound 5}: The most we can assert with proof is that for any $\epsilon > 0$ we have
\[
\mbox{\as}:\ \tilde{\gamma}_m
\leq (1 + o(1)) \tfrac{1}{6} \tilde{\rho}_m^3
= O\left( m^{\frac{9}{4} + \epsilon} \right),
\]
where the inequality follows from \refT{T:ub} and the big-oh bound from two results in~\cite{Fillboundary(2018)}:\ the big-oh bound (4.7) and Prop.\ 5.1(b1).  But for $d = 2, 3$ the situation improves relative to~\eqref{as bound 5}.  For $d = 3$, for any $\epsilon > 0$ we know
\[
\mbox{\as}:\ \tilde{\gamma}_m = 2 \tilde{\rho}_m + 1 = O\left( m^{\frac{3}{4} + \epsilon} \right),
\]
where the equality follows from \refC{C:3d} and the big-oh bound from the same two results in~\cite{Fillboundary(2018)} as for $d = 4$.  For $d = 2$, we are able to prove a big-oh bound that matches~\eqref{conjecture} in magnitude:
\[
\mbox{\as}:\ \tilde{\gamma}_m = \tilde{\rho}_m + 1 = O(m^{1 / 2}),
\]
where the equality follows from \refC{C:3d} and the big-oh bound again from two results in~\cite{Fillboundary(2018)}:\ the end of Rmk.\ 4.3 and Prop.\ 5.1(b1).
In fact, the proof discussed at the end of \cite[Rmk.\ 4.3]{Fillboundary(2018)} gives the sharper result
\[
\mbox{\as}:\ \tilde{\gamma}_m = \tilde{\rho}_m + 1 \leq (1 + o(1)) e \L T_m \leq (1 + o(1)) \sqrt{2} e m^{1/2};
\]
the multiplicative constant here is larger than the conjectured constant in~\eqref{conjecture} only by a factor of~$e$.
However, we
concede that the result for $d = 2$ might have less use in analyzing the cost of generating records than for $d \geq 3$, since in that case one has a simpler representation of the record-setting region that obviates use of the importance-sampling algorithm (recall \refR{R:bivariate}).
\end{remark}

\begin{acks}
We thank Vince Lyzinski and Fred Torcaso for helpful comments.
\end{acks}

\bibliography{records}
\bibliographystyle{plain}

\end{document}

%% file: figure0.tex
\begin{adjustbox}{max totalsize={.9\textwidth}{.7\textheight},center}
\begin{tikzpicture}[scale=7]
%
%
\draw[-{Triangle},thick,color=black] (0,0)--(0,1.2);
\draw[-{Triangle},thick,color=black] (0,0)--(1.2,0);

%
%
\draw[thick,color=black](0,.76)--
(.33,.76)--(.33,.49) --
(.51,.49)--(.51,.28) --
(.62,.28)--(.62,.11) --
(.86,.11)-- (.86,.0)
;
%
%
\draw[thick,color=black](0,1.)--(1.,1.)--(1.,0.);

%
%
\filldraw [black]
(0.0,.76) circle (.25pt)
(.33,.49) circle (.25pt)
(.51,.28) circle (.25pt)
(.62,.11) circle (.25pt)
(.86,0.) circle (.25pt)
;
%
%
\draw [black]
(.33,.76) circle (.25pt)
(.51,.49) circle (.25pt)
(.62,.28) circle (.25pt)
(.86,.11) circle (.25pt)
;
%
%
%
%
%
%
\def\mx{.04}
\def\my{.04}
\draw (0.-\mx,.76-\my) node[color=black] {\footnotesize $g^{(1)}$};
\draw (.33-\mx,.49-\my) node[color=black] {\footnotesize $g^{(2)}$};
\draw (.51-\mx,.28-\my) node[color=black] {\footnotesize  $g^{(3)}$};
\draw (.62-\mx,.11-\my) node[color=black] {\footnotesize $g^{(4)}$};
\draw (.86-\mx,0.-\my) node[color=black] {\footnotesize $g^{(5)}$};

%
%
\def\mmx{.05}
\def\mmy{.04}
\draw (.33+\mmx,.76+\mmy) node[color=black] {\footnotesize $r^{(1)}$};
\draw (.51+\mmx,.49+\mmy) node[color=black] {\footnotesize $r^{(2)}$};
\draw (.62+\mmx,.28+\mmy) node[color=black] {\footnotesize $r^{(3)}$};
\draw (.86+\mmx,.11+\mmy) node[color=black] {\footnotesize $r^{(4)}$};

%
%
\draw (.75,.64) node[color=black] {\Large $S$};

\draw(1.05,1.05) node[color=black]{\footnotesize $(1,1)$};

%
%

%
%
\draw (1.20,-.05) node[color=black] {\footnotesize $x_1$};
\draw (-.05,1.20) node[color=black] {\footnotesize $x_2$};

%
%
\draw[dashed,color=black](.51,.28)->(.51,1.00);

\draw[dashed,color=black](.33,.49)->(.33,1.00);

\draw[dashed,color=black](.62,.11)->(.62,1.00);


\draw[dashed,color=black](.86,.11)->(.86,1.);

\end{tikzpicture}
\end{adjustbox}

%% file: figure1.tex
\begin{adjustbox}{max totalsize={.9\textwidth}{.7\textheight},center}
\begin{tikzpicture}[scale=7]
\draw[color=black,solid,thick](-0.2575755771623076,-1.9142739173960515)--(0.5362581085663101,-1.455953825248841);
\draw[color=black,solid,thick](0.5362581085663101,-1.455953825248841)--(-0.43444147459620586,-0.8955201596077051);
\draw[color=black,solid,thick](0.5362581085663101,-1.455953825248841)--(0.5362581085663101,-1.636381118407808);
\draw[color=black,solid,thick](-0.2575755771623076,-1.9142739173960515)--(-1.942852250259816,-0.941278976524208);
\draw[color=black,solid,thick](-1.942852250259816,-0.941278976524208)--(-1.3624390450436483,-0.60617725625144);
\draw[color=black,solid,thick](-1.942852250259816,-0.941278976524208)--(-1.942852250259816,-1.121706269683175);
\draw[color=black,solid,thick](-0.2575755771623076,-1.9142739173960515)--(-0.2575755771623076,-2.0947012105550185);
\draw[color=black,solid,thick](-0.2575755771623076,-2.0947012105550185)--(0.5362581085663101,-1.636381118407808);
\draw[color=black,solid,thick](-0.2575755771623076,-2.0947012105550185)--(-1.942852250259816,-1.121706269683175);
\draw[color=black,solid,thick](1.1865525987720142,0.06663339463679119)--(1.6313637441374893,0.32344522915209406);
\draw[color=black,solid,thick](1.6313637441374893,0.32344522915209406)--(1.4227264509264796,0.4439020265504653);
\draw[color=black,solid,thick](1.6313637441374893,0.32344522915209406)--(1.6313637441374893,-0.4304951349833537);
\draw[color=black,solid,thick](1.1865525987720142,0.06663339463679119)--(0.2513166656967584,0.6065921123536011);
\draw[color=black,solid,thick](0.2513166656967584,0.6065921123536011)--(0.3117430076602592,0.6414792771523726);
\draw[color=black,solid,thick](0.2513166656967584,0.6065921123536011)--(0.2513166656967584,0.36910122321068184);
\draw[color=black,solid,thick](1.1865525987720142,0.06663339463679119)--(1.1865525987720142,-0.6873069694986566);
\draw[color=black,solid,thick](1.1865525987720142,-0.6873069694986566)--(1.6313637441374893,-0.4304951349833537);
\draw[color=black,solid,thick](1.1865525987720142,-0.6873069694986566)--(0.9978411846641773,-0.5783543837643418);
\draw[color=black,solid,thick](-0.10188361164666793,2.115297774717315)--(0.6961278110622334,2.5760298844346954);
\draw[color=black,solid,thick](0.6961278110622334,2.5760298844346954)--(0.0,2.977939463541861);
\draw[color=black,solid,thick](0.6961278110622334,2.5760298844346954)--(0.6961278110622334,1.9561916149441554);
\draw[color=black,solid,thick](-0.10188361164666793,2.115297774717315)--(-0.7980114227089014,2.5172073538244804);
\draw[color=black,solid,thick](-0.7980114227089014,2.5172073538244804)--(0.0,2.977939463541861);
\draw[color=black,solid,thick](-0.7980114227089014,2.5172073538244804)--(-0.7980114227089014,0.5670905271157696);
\draw[color=black,solid,thick](-0.10188361164666793,2.115297774717315)--(-0.10188361164666793,0.16518094800860406);
\draw[color=black,solid,thick](-0.10188361164666793,0.16518094800860406)--(0.2513166656967584,0.36910122321068184);
\draw[color=black,solid,thick](-0.10188361164666793,0.16518094800860406)--(-0.7980114227089014,0.5670905271157696);
\draw[color=black,solid,thick](1.0550507459120508,1.4061056218090915)--(1.3519777228080998,1.577536491849687);
\draw[color=black,solid,thick](1.3519777228080998,1.577536491849687)--(0.6961278110622334,1.9561916149441554);
\draw[color=black,solid,thick](1.3519777228080998,1.577536491849687)--(1.3519777228080998,0.6654418574334658);
\draw[color=black,solid,thick](1.0550507459120508,1.4061056218090915)--(0.39920083416618446,1.7847607449035596);
\draw[color=black,solid,thick](0.39920083416618446,1.7847607449035596)--(0.6961278110622334,1.9561916149441554);
\draw[color=black,solid,thick](0.39920083416618446,1.7847607449035596)--(0.39920083416618446,0.8726661104873383);
\draw[color=black,solid,thick](1.0550507459120508,1.4061056218090915)--(1.0550507459120508,0.4940109873928702);
\draw[color=black,solid,thick](1.0550507459120508,0.4940109873928702)--(1.3519777228080998,0.6654418574334658);
\draw[color=black,solid,thick](1.0550507459120508,0.4940109873928702)--(0.39920083416618446,0.8726661104873383);
\draw[color=black,solid,thick](-0.6478619551086559,-0.242005219307163)--(-0.201089630713166,0.01593890244904783);
\draw[color=black,solid,thick](-0.201089630713166,0.01593890244904783)--(-0.9156667206481585,0.4285001776797551);
\draw[color=black,solid,thick](-0.201089630713166,0.01593890244904783)--(-0.201089630713166,-0.42984854104428194);
\draw[color=black,solid,thick](-0.6478619551086559,-0.242005219307163)--(-1.3624390450436483,0.17055605592354428);
\draw[color=black,solid,thick](-1.3624390450436483,0.17055605592354428)--(-0.9156667206481585,0.4285001776797551);
\draw[color=black,solid,thick](-1.3624390450436483,0.17055605592354428)--(-1.3624390450436483,-0.60617725625144);
\draw[color=black,solid,thick](-0.6478619551086559,-0.242005219307163)--(-0.6478619551086559,-1.0187385314821473);
\draw[color=black,solid,thick](-0.6478619551086559,-1.0187385314821473)--(-0.43444147459620586,-0.8955201596077051);
\draw[color=black,solid,thick](-0.6478619551086559,-1.0187385314821473)--(-1.3624390450436483,-0.60617725625144);
\draw[color=black,solid,thick](1.190670477917273,-1.502833114172986)--(2.3396890424484713,-0.8394469366370213);
\draw[color=black,solid,thick](2.3396890424484713,-0.8394469366370213)--(1.6313637441374893,-0.4304951349833537);
\draw[color=black,solid,thick](2.3396890424484713,-0.8394469366370213)--(2.3396890424484713,-1.3508200984776428);
\draw[color=black,solid,thick](1.190670477917273,-1.502833114172986)--(-0.43444147459620586,-0.5645742909260506);
\draw[color=black,solid,thick](-0.43444147459620586,-0.5645742909260506)--(-0.201089630713166,-0.42984854104428194);
\draw[color=black,solid,thick](-0.43444147459620586,-0.5645742909260506)--(-0.43444147459620586,-0.8955201596077051);
\draw[color=black,solid,thick](1.190670477917273,-1.502833114172986)--(1.190670477917273,-2.0142062760136077);
\draw[color=black,solid,thick](1.190670477917273,-2.0142062760136077)--(2.3396890424484713,-1.3508200984776428);
\draw[color=black,solid,thick](1.190670477917273,-2.0142062760136077)--(0.5362581085663101,-1.636381118407808);
\draw[color=black,solid,thick](1.0383416475245055,0.40267039049296427)--(1.4227264509264796,0.6245950602094955);
\draw[color=black,solid,thick](1.4227264509264796,0.6245950602094955)--(1.3519777228080998,0.6654418574334658);
\draw[color=black,solid,thick](1.4227264509264796,0.6245950602094955)--(1.4227264509264796,0.4439020265504653);
\draw[color=black,solid,thick](1.0383416475245055,0.40267039049296427)--(0.3117430076602592,0.8221723108114029);
\draw[color=black,solid,thick](0.3117430076602592,0.8221723108114029)--(0.39920083416618446,0.8726661104873383);
\draw[color=black,solid,thick](0.3117430076602592,0.8221723108114029)--(0.3117430076602592,0.6414792771523726);
\draw[color=black,solid,thick](1.0383416475245055,0.40267039049296427)--(1.0383416475245055,0.22197735683393405);
\draw[color=black,solid,thick](1.0383416475245055,0.22197735683393405)--(1.4227264509264796,0.4439020265504653);
\draw[color=black,solid,thick](1.0383416475245055,0.22197735683393405)--(0.3117430076602592,0.6414792771523726);
\draw[color=black,solid,thick](0.5269856093814939,-0.33375350191070674)--(0.9978411846641773,-0.061904908771813316);
\draw[color=black,solid,thick](0.9978411846641773,-0.061904908771813316)--(0.2513166656967584,0.36910122321068184);
\draw[color=black,solid,thick](0.9978411846641773,-0.061904908771813316)--(0.9978411846641773,-0.5783543837643418);
\draw[color=black,solid,thick](0.5269856093814939,-0.33375350191070674)--(-0.9156667206481585,0.4991622091789538);
\draw[color=black,solid,thick](-0.9156667206481585,0.4991622091789538)--(-0.7980114227089014,0.5670905271157696);
\draw[color=black,solid,thick](-0.9156667206481585,0.4991622091789538)--(-0.9156667206481585,0.4285001776797551);
\draw[color=black,solid,thick](0.5269856093814939,-0.33375350191070674)--(0.5269856093814939,-0.8502029769032352);
\draw[color=black,solid,thick](0.5269856093814939,-0.8502029769032352)--(0.9978411846641773,-0.5783543837643418);
\draw[color=black,solid,thick](0.5269856093814939,-0.8502029769032352)--(-0.201089630713166,-0.42984854104428194);
\draw[->,color=black,solid,thick](-1.942852250259816,-1.121706269683175)--(-2.525707925337761,-1.4582181505881275);
\draw[->,color=black,solid,thick](2.3396890424484713,-1.3508200984776428)--(3.041595755183013,-1.7560661280209358);
\draw[->,color=black,solid,thick](0.0,2.977939463541861)--(0.0,3.8713213026044193);
\filldraw [red]
(-0.2575755771623076,-1.9142739173960515) circle (1.25pt)
(1.1865525987720142,0.06663339463679119) circle (1.25pt)
(-0.10188361164666793,2.115297774717315) circle (1.25pt)
(1.0550507459120508,1.4061056218090915) circle (1.25pt)
(-0.6478619551086559,-0.242005219307163) circle (1.25pt)
(1.190670477917273,-1.502833114172986) circle (1.25pt)
(1.0383416475245055,0.40267039049296427) circle (1.25pt)
(0.5269856093814939,-0.33375350191070674) circle (1.25pt)
;\filldraw [green]
(-0.43444147459620586,-0.8955201596077051) circle (1.25pt)
(0.3117430076602592,0.6414792771523726) circle (1.25pt)
(0.2513166656967584,0.36910122321068184) circle (1.25pt)
(0.9978411846641773,-0.5783543837643418) circle (1.25pt)
(0.39920083416618446,0.8726661104873383) circle (1.25pt)
(-0.201089630713166,-0.42984854104428194) circle (1.25pt)
;\filldraw [blue]
(0.6961278110622334,1.9561916149441554) circle (1.25pt)
(1.3519777228080998,0.6654418574334658) circle (1.25pt)
(1.4227264509264796,0.4439020265504653) circle (1.25pt)
(1.6313637441374893,-0.4304951349833537) circle (1.25pt)
(-0.7980114227089014,0.5670905271157696) circle (1.25pt)
(-0.9156667206481585,0.4285001776797551) circle (1.25pt)
(-1.3624390450436483,-0.60617725625144) circle (1.25pt)
(0.5362581085663101,-1.636381118407808) circle (1.25pt)
;\filldraw [violet]
(-1.942852250259816,-1.121706269683175) circle (1.25pt)
(2.3396890424484713,-1.3508200984776428) circle (1.25pt)
(0.0,2.977939463541861) circle (1.25pt)
;\draw (-2.4,-1.6) node[color=black, scale=3.5]{ $x_1$};
\draw (2.9,-1.9) node[color=black, scale=3.5]{ $x_2$};
\draw (0.2,3.5) node[color=black, scale=3.4]{ $x_3$};
\draw[->,color=black,dashed,ultra thick](-0.201089630713166,-0.42984854104428194)--(-1.933140438282043,-1.429848541044282);
\draw[->,color=black,dashed,ultra thick](-0.201089630713166,-0.42984854104428194)--(1.5309611768557114,-1.429848541044282);
\draw[->,color=black,dashed,ultra thick](-0.201089630713166,-0.42984854104428194)--(-0.201089630713166,1.5701514589557182);
%
%
\draw (-.2+.10,-.43+.08) node[color=black, scale=3.4]{ $g$};

\end{tikzpicture}
\end{adjustbox}

%% file: figure2.tex
\begin{adjustbox}{max totalsize={.6\textwidth}{.6\textheight},center}
\begin{tikzpicture}[scale=7]
\draw[color=black,solid,thick](-0.5196152422706635,-1.4999999999999998)--(2.1650635094610964,0.050000000000000044);
\draw[color=black,solid,thick](2.1650635094610964,0.050000000000000044)--(1.3856406460551018,0.5);
\draw[color=black,solid,thick](2.1650635094610964,0.050000000000000044)--(2.1650635094610964,-1.25);
\draw[color=black,solid,thick](-0.5196152422706635,-1.4999999999999998)--(-2.68467875173176,-0.25);
\draw[color=black,solid,thick](-2.68467875173176,-0.25)--(-1.9918584287042087,0.15000000000000013);
\draw[color=black,solid,thick](-2.68467875173176,-0.25)--(-2.68467875173176,-1.55);
\draw[color=black,solid,thick](-0.5196152422706635,-1.4999999999999998)--(-0.5196152422706635,-2.8);
\draw[color=black,solid,thick](-0.5196152422706635,-2.8)--(2.1650635094610964,-1.25);
\draw[color=black,solid,thick](-0.5196152422706635,-2.8)--(-2.68467875173176,-1.55);
\draw[color=black,solid,thick](-0.6062177826491069,-0.19999999999999996)--(1.3856406460551018,0.95);
\draw[color=black,solid,thick](1.3856406460551018,0.95)--(0.6928203230275509,1.35);
\draw[color=black,solid,thick](1.3856406460551018,0.95)--(1.3856406460551018,0.5);
\draw[color=black,solid,thick](-0.6062177826491069,-0.19999999999999996)--(-1.9918584287042087,0.6000000000000001);
\draw[color=black,solid,thick](-1.9918584287042087,0.6000000000000001)--(-0.9526279441628825,1.2);
\draw[color=black,solid,thick](-1.9918584287042087,0.6000000000000001)--(-1.9918584287042087,0.15000000000000013);
\draw[color=black,solid,thick](-0.6062177826491069,-0.19999999999999996)--(-0.6062177826491069,-0.6499999999999999);
\draw[color=black,solid,thick](-0.6062177826491069,-0.6499999999999999)--(1.3856406460551018,0.5);
\draw[color=black,solid,thick](-0.6062177826491069,-0.6499999999999999)--(-1.9918584287042087,0.15000000000000013);
\draw[color=black,solid,thick](-0.2598076211353316,1.4499999999999997)--(0.6928203230275509,2.0);
\draw[color=black,solid,thick](0.6928203230275509,2.0)--(0.3377499074759311,2.205);
\draw[color=black,solid,thick](0.6928203230275509,2.0)--(0.6928203230275509,1.35);
\draw[color=black,solid,thick](-0.2598076211353316,1.4499999999999997)--(-0.9526279441628825,1.8499999999999999);
\draw[color=black,solid,thick](-0.9526279441628825,1.8499999999999999)--(-0.3897114317029974,2.175);
\draw[color=black,solid,thick](-0.9526279441628825,1.8499999999999999)--(-0.9526279441628825,1.2);
\draw[color=black,solid,thick](-0.2598076211353316,1.4499999999999997)--(-0.2598076211353316,0.7999999999999999);
\draw[color=black,solid,thick](-0.2598076211353316,0.7999999999999999)--(0.6928203230275509,1.35);
\draw[color=black,solid,thick](-0.2598076211353316,0.7999999999999999)--(-0.9526279441628825,1.2);
\draw[color=black,solid,thick](-0.05196152422706629,2.2800000000000002)--(0.3377499074759311,2.5050000000000003);
\draw[color=black,solid,thick](0.3377499074759311,2.5050000000000003)--(0.0,2.7);
\draw[color=black,solid,thick](0.3377499074759311,2.5050000000000003)--(0.3377499074759311,2.205);
\draw[color=black,solid,thick](-0.05196152422706629,2.2800000000000002)--(-0.3897114317029974,2.475);
\draw[color=black,solid,thick](-0.3897114317029974,2.475)--(0.0,2.7);
\draw[color=black,solid,thick](-0.3897114317029974,2.475)--(-0.3897114317029974,2.175);
\draw[color=black,solid,thick](-0.05196152422706629,2.2800000000000002)--(-0.05196152422706629,1.98);
\draw[color=black,solid,thick](-0.05196152422706629,1.98)--(0.3377499074759311,2.205);
\draw[color=black,solid,thick](-0.05196152422706629,1.98)--(-0.3897114317029974,2.175);
\draw[->,color=black,solid,thick](-2.68467875173176,-1.55)--(-3.4900823772512877,-2.015);
\draw[->,color=black,solid,thick](2.1650635094610964,-1.25)--(2.8145825622994254,-1.625);
\draw[->,color=black,solid,thick](0.0,2.7)--(0.0,3.5100000000000002);
\filldraw [red]
(-0.5196152422706635,-1.4999999999999998) circle (1.75pt)
(-0.6062177826491069,-0.19999999999999996) circle (1.75pt)
(-0.2598076211353316,1.4499999999999997) circle (1.75pt)
(-0.05196152422706629,2.2800000000000002) circle (1.75pt)
;\filldraw [green]
;\filldraw [blue]
(0.3377499074759311,2.205) circle (1.75pt)
(0.6928203230275509,1.35) circle (1.75pt)
(1.3856406460551018,0.5) circle (1.75pt)
(-0.3897114317029974,2.175) circle (1.75pt)
(-0.9526279441628825,1.2) circle (1.75pt)
(-1.9918584287042087,0.15000000000000013) circle (1.75pt)
;\filldraw [violet]
(-2.68467875173176,-1.55) circle (1.75pt)
(2.1650635094610964,-1.25) circle (1.75pt)
(0.0,2.7) circle (1.75pt)
;\draw (-3.3,-2.2) node[color=black, scale=5]{ $x_1$};
\draw (2.7,-1.8) node[color=black, scale=5]{ $x_2$};
\draw (0.2,3.5) node[color=black, scale=5]{ $x_3$};

\end{tikzpicture}
\end{adjustbox}